%% file: Main.tex
\documentclass[11pt]{article}

\makeatletter

\renewcommand\section{\@startsection {section}{1}{\z@}
                                       {-3.5ex \@plus -1ex \@minus -.2ex}
                                       {2.3ex \@plus.2ex}
                                       {\normalfont\fontfamily{phv}\fontsize{16}{19}\bfseries}}
\renewcommand\subsection{\@startsection{subsection}{2}{\z@}
                                         {-3.25ex\@plus -1ex \@minus -.2ex}
                                         {1.5ex \@plus .2ex}
                                         {\normalfont\fontfamily{phv}\fontsize{14}{17}\bfseries}}
\renewcommand\subsubsection{\@startsection{subsubsection}{3}{\z@}
                                        {-3.25ex\@plus -1ex \@minus -.2ex}
                                         {1.5ex \@plus .2ex}
                                         {\normalfont\normalsize\fontfamily{phv}\fontsize{14}{17}\selectfont}}
\makeatother
   
\usepackage{amsmath}
\usepackage{graphicx}
\usepackage{enumerate}
\usepackage{xcolor}
\usepackage{natbib}
\usepackage{url}
\usepackage{amsfonts, amsthm, latexsym, amssymb}
\usepackage{mathtools}
\usepackage{lineno}

\usepackage{algorithmic} 
\usepackage{comment}
\usepackage[ruled]{algorithm2e}
\SetAlgoNlRelativeSize{-1}
\usepackage{enumitem}
\usepackage{multirow}
\usepackage{mathrsfs}
\usepackage{arydshln}
\usepackage[resetlabels]{multibib}
\newcites{Supp}{References}
\usepackage{natbib} 
 	
\begin{document}

\def\spacingset#1{\renewcommand{\baselinestretch}{#1}\small\normalsize} \spacingset{1}
		
\title{\bf When Shortages Lead to Export Restrictions: \\
            A Computational Study}
		\author{Martha L. Sabogal De La Pava $^a$, Emily L. Tucker $^{a}$, and Yongjia Song $^{a}$ \\
			$^a$ Industrial Engineering Department, Clemson University, Clemson, United States \\ 
          }
\date{}
\maketitle

\begin{abstract}
Globalization has enabled companies to leverage cost efficiencies; however, it has increased exposure to disruption risks that threaten supply stability. Among these, export bans have emerged as a systemic challenge, often arising as a secondary effect of conventional supply capacity disruptions. The pharmaceutical industry is particularly vulnerable to bans because of the criticality of keeping life-saving medications available, yet such disruptions also jeopardize patient health globally as well as challenge company operations. This paper proposes a supply chain design model that incorporates the relationship between conventional disruption risks (quality-related capacity failures and natural disasters) and export bans that may be induced when the drug is projected to be in short supply. 
The model is a two-stage stochastic mixed integer program, where the systemic setting yields a second-stage problem with one binary variable and continuous variables. The structure motivates the study of three tailored solution methods: the Alternating Integer L-shaped Cut, Alternating Disjunctive Cut, and Alternating Bilinear Cut methods. In the nearly-continuous second stage context, the Alternating Integer L-shaped cuts perform best. Using an oncology drug case study, we study effects of systemic, induced risk and quantify the value of incorporating disruptions into supply chain design.
\end{abstract}
			
\spacingset{1.5}

\section{Introduction} \label{s:intro}
Pharmaceutical supply chains are global operations that serve patients worldwide. The generic pharmaceutical industry, where multiple companies may produce equivalent products, operates in a particularly lean setting due to the underlying economics, where redundancy may not be cost-effective \citep{Zhao-tutorial, Tucker2020}. The supply chains also tend to be more rigid than in other industries because of safety regulations. As a result, the ability to meet consumer demands is highly dependent on the structure of the supply chain, which is decided far ahead of disruptions.

The global setting exposes generic pharmaceutical supply chains to conventional disruptions such as natural disasters or quality disruptions, where plants may be closed for years due to contamination or other risks. For example, contamination risks led to plant closures in India and widespread oncology drug shortages \citep{WSJ-Mosbergen-Quality2023}, and a tornado shut down a major plant in North Carolina, wiping products off of wholesaler shelves \citep{NBC-Lovelace-Tornado2023}. These shortages may last for years because of the challenge of recovering and qualifying facilities.

Yet, the criticality of life-saving medications has also led to countries to impose or threaten to impose export bans during shortages. These are intended to keep the limited supply available in-country. Bulgarian regulations, for example, stipulate that if the availability of a specific drug falls below 65\% of the domestic demand for a period of one month, local authorities may impose export bans on that drug \citep{BulgariaDrugAgency2024}; they imposed export bans on insulin and children's antibiotics in 2024 \citep{Euractiv-Nokolov-Export2024}. India, Canada, and Italy have threatened similar policies. The risk of export bans imposed by multiple countries due to shortages behaves as an \textit{induced} systemic disruption, where a conventional disruption, such as a hurricane, may trigger a large-scale, global disruption. 

Existing supply chain models that include disruptions generally parameterize a single, generic disruption parameter \citep{Snyder2016,Aldrighetti2021} and fail to account for simultaneous disruption risks \citep{Ghadge2022}. In the current global setting, supply chain risk is multi-faceted \citep{Shekarian2021}. The uniqueness of the pharmaceutical industry, including the potential for induced export bans, motivates a new study on global supply chain planning with disaggregated disruption dynamics, including the new prominence of induced disruptions. These disruptions include regional risks (e.g., natural disasters); facility-specific risks that are affected by quality investments; and geopolitical risks such as export bans that are affected by supply decisions and projected shortages. Disaggregating disruptions by source of strain may be critical to develop appropriate supply chain plans in the new era of integrated global disruption risks.

In this paper, we propose a two-stage stochastic mixed-integer programming model for the problem of global pharmaceutical supply chain design under conventional and systemic disruptions. It includes: (i) quality-related capacity disruptions, (ii) correlated disruptions due to natural disasters, and (iii) export bans that are induced by projected drug shortages.

A key feature of the proposed model is its approach to mathematically formulate the systemic shock of the export ban activation, which results in a special second-stage problem containing continuous variables and a single binary variable per scenario. The general case of two-stage stochastic programming with a second-stage mixed-integer program is challenging, and existing methods have been studied for problems with a general number of binary and/or integer variables \citep{van-der-2024}. The unique structure of the pharmaceutical export ban setting, however, allows us to develop tailored cutting plane generation methods. In this work, we evaluate the computational performance of three cutting plane methods based on a Benders-style decomposition framework, and the key contributions of this paper are as follows.
\begin{itemize}[noitemsep,topsep=0pt]
    \item It introduces the first supply chain design model in which export bans are a response to projected shortages, stemming from facility quality failures, natural disasters, and supply-demand mismatches.
    \item It evaluates the computational performance of three cutting plane methods for solving a two-stage mixed-integer stochastic program with pure binary first-stage decisions and one binary variable and continuous variables per scenario in the second stage. 
    \item It assesses the value of incorporating systemic and conventional disruption risks into the supply chain planning context by quantifying the benefits of planning for individual disruption risks and their combinations. It further identifies the types of disruptions that have the most influence on supply chain design decisions using a case study in the pharmaceutical industry.
\end{itemize}

\section{Literature Review} \label{s:LR}
We first review literature related to risks of geopolitical strain, facility quality failures, and natural disasters in the context of supply chain design. Then, we present literature on algorithms for solving two-stage stochastic mixed-integer programs.

Research on supply chain design under disruption risks is extensive (see, e.g., reviews by~\cite{Katsaliaki2022,Suryawanshi2022}). In this context, disruptions are commonly categorized based on geographic impact as \textit{local} (e.g., failure in a facility), \textit{regional} (e.g., natural disasters affecting multiple components in a supply chain), and \textit{global} (e.g., major strikes and pandemics) \citep{Roshani2024}. The majority of existing studies focus on local disruptions \citep{Roshani2024}. Typically, design optimization models do not integrate multiple types of disruptions \citep{Ghadge2022} and are often included using a single parameter \citep{Snyder2016,Aldrighetti2021}, representing partial or complete failures of suppliers, facilities, links, and networks \citep{HARTNIBBRIG2025}.

Geopolitical risks such as wars, political instability, and terrorism may lead to strikes, trade restrictions, and regulatory barriers that may disrupt supply chain operations \citep{Roscoe2022}. Despite their growing relevance, research on these risks remains limited \citep{Bednarski2025}. Most studies are in the operations management literature, including frameworks such as diffusion style models, variational inequality approaches, and empirical studies \citep{Bednarski2025}; and focus on specific industries, e.g., lithium-ion batteries \citep{REN2024}, medical goods \citep{Grassia2021}, semiconductors \citep{TSE2024}, and agriculture \citep{Hassani2025}. Optimization-based models that incorporate geopolitical risks into supply chain design are largely underexplored. One exception is proposed by \cite{SABOGAL2025}. They model export bans as an uncertain, exogenous parameter based on the realized availability of raw materials; this represents one leading indicator of potential shortages. Building on this idea, our work increases realism by modeling export bans as triggered when shortages themselves are projected to occur, which reflects what decision-makers may focus on when making decisions to activate bans. While shortages are regularly included within optimization models with a penalty cost \citep{GUO2025}, our approach links projected shortages to a new risk: induced disruptions.

Quality is a traditional performance indicator for supply chain design \citep{Maharjan2022}, and extensive research has focused on integrating quality concerns into supply chain design models (e.g., \cite{FRANCA2010,RAMEZANI2013}). In the context of facility disruptions, some studies have explicitly considered quality failures and quality levels in their modeling approaches (e.g., \cite{MADADI2014,Tucker2022}). These studies show the value of considering quality-related factors in enhancing supply chain reliability and reducing shortages. Our paper includes facility disruptions as a function of the selected quality levels, where higher-quality facilities are assumed to be more reliable. This relates to the reliable network design and infrastructure protection literature, in which a facility's post-disruption capacity depends on protection investment decisions (e.g., \cite{Bhuiyan2020,HU2025}). One difference, however, between infrastructure and supply chain planning is that infrastructure planning focuses on access, while for-profit supply chains (like our problem setting) prioritize profitability.

Natural disasters can induce correlations in disruptions across facilities \citep{Wang2024}, and they are often studied in the context of reliable facility location. The two main approaches to incorporate correlated disruptions are the use of supporting station models and shared hazards with propagation effects. In a supporting station modeling approach, the interdependent and positively correlated facility failures are transformed into a supporting structure with i.i.d. disruptions. A pioneering paper in this line of research is from \cite{Li-Fan2013} with extensions that consider probabilistic \citep{Xie-Li2015} and generalized correlation patterns \citep{Xie-An2019}. A shared hazard exposure and propagation approach is proposed by \cite{Li-Ouyang2010}, considering spatially correlated disruptions (e.g., based on distance) with site-dependent disruption probabilities in a continuum approximation model. \cite{Lu-Ran2015} and \cite{Li2018} extend this research by considering a distributionally robust optimization approach and a stochastic programming framework, respectively. Alternative approaches are also possible, e.g., \cite{Liberatore2012} present a tri-level fortification model with a correlation matrix. Research on pharmaceutical supply chains with correlated disruptions is limited \citep{Hasani2016,Pariazar2017}, and in general, the challenge of including correlated disruptions in optimization models is computational tractability \citep{Xie-An2019}. These aspects motivate further study on computational performance when correlated disruptions are integrated into new model structures, as in this paper, where strategic facility choices with correlated disruption risks may induce a systemic disruption.

The model we present is a two-stage stochastic mixed-integer program, and this structure has been used to model a wide range of problems, including the supply chain design problem \citep{KSen-tutorial}. These models can be challenging to solve as they combine the complexities of solving integer programs, such as non-convexity and discontinuity, with the large-scale context of stochastic programs \citep{Keller2012}. 

When integer variables are only in the first stage and the second-stage problems are linear programs, Benders decomposition (also called the L-shaped method) can be applied since the recourse function is convex and can be approximated by piecewise linear functions \citep{Birge_Louveaux_2011}. When there are integer variables in both stages, these properties do not hold, and to solve these programs, modified Benders decomposition and dual decomposition are the most common approaches \citep{Chen2022}.

There is growing research on cutting planes to improve Benders-style decompositions. We focus our discussion on these methods. A seminal work by \cite{Laporte1993} proposes the integer L-shaped method for stochastic programs with binary first-stage variables and arbitrary second-stage variables. When the second-stage problem (also called subproblem) is a mixed-integer program, the integer L-shaped method requires the exact solution of the subproblem. To improve this method, \cite{Angulo2016} propose the alternating cuts strategy and the use of a cut-generating linear program-based optimality cuts. The alternating strategy aims to reduce the expensive evaluation of the second-stage value function by alternating between solving linear and mixed-integer subproblems. The new optimality cuts enhance the approximation of the value function. 

Another approach is disjunctive decomposition, where disjunctive programming theory is used to convexify the mixed-integer feasible region of the second-stage problem \citep{Sen-Higle2005}. Disjunctive decomposition solves the linear relaxation of the subproblem. This provides a computational advantage over the integer L-shaped method, which requires solving the subproblem as a mixed-integer program. However, when solving the subproblem relaxation, if the solution does not satisfy the integrality constraints, a cut-generating problem needs to be solved to find a disjunctive cut \citep{Keller2012}. In particular, \cite{Sen-Sherali2006} propose to combine disjunctive decomposition with the branch-and-cut method; the method adds disjunctive cuts formed on the subproblems' value functions to the master problem. Gomory cuts can also be used within a Benders framework. For instance, \cite{Gade2014} propose an algorithm to solve stochastic programs with pure binary first-stage and integer second-stage variables, and \cite{Zhang-K2014} propose an algorithm for the general pure integer case in both stages. \cite{van-der-2024} focus on a general two-stage stochastic mixed-integer model and introduce scaled cuts that can be used to recover the convex envelope of the expected second-stage objective function. 

In summary, this paper contributes to the literature with a supply chain design model that considers multiple sources of disruptions within a unified framework, including emerging threats of shortage-induced export bans, and provides a computational study of three cutting plane methods to solve a two-stage stochastic program with pure binary first-stage variables and a second stage with continuous variables and one binary variable per scenario.

\section{Model Description} \label{s:model}

The model takes the perspective of a generic pharmaceutical company that seeks to maximize its expected profit by producing and selling a drug under conventional and systemic disruptions and operational uncertainty. It is formulated as a two-stage stochastic program. In the first stage, the company selects the location and quality level of its manufacturing plants ($y_{jl} \in \{0,1\}$ for each candidate plant $j\in J$ and quality level $l\in L$). Then, all uncertainty is realized. In each second-stage scenario ($\omega \in \Omega)$, the company orders raw materials ($u_{ij}^\omega$) based on the production plan ($x_j^\omega$), transships between plants ($\tau_{jj^\prime}^\omega$), and produces the finished product and ships to customers ($v_{jk}^\omega$), indexed by supplier $i \in I$, candidate plant $j\in J$, and country with demand $k\in K$. 

The new feature is the risk of induced export bans. The company's ability to ship finished product from a plant where it was manufactured to other countries (i.e., exports) depends on if the country in which it was manufactured allows exports of the drug. We model export ban activation following the threshold-based design of drug shortage-related export ban policies (e.g., \cite{BulgariaDrugAgency2024}). At a high-level, if a country $j$ has an exogenous export ban policy ($j \in \bar{J} \subseteq J$), then the country $j$ will impose an export ban if the \textit{projected global drug supply} is below a threshold on the \textit{nominal global drug availability}.  The projected global drug supply is based on the (i) company's projected production, considering conventional disruptions and operational uncertainty and (ii) the exogenous drug production from other companies that produce the drug. Nominal global drug availability is based on a model without disruptions. Because export ban thresholds may be written into law, the activation is modeled as a second-stage decision; the company whose perspective we take has information on how their choices will contribute to the projected global supply and whether there will be a resulting systemic disruption. Export allowability is expressed with the decision variable, $z^\omega \in \{0,1\}$ where a value of 1 allows exports and 0 bans exports from all countries in set $\bar{J}$ in scenario $\omega \in \Omega$, representing a systemic, induced disruption. The resulting model is a two-stage mixed-integer stochastic program with binary first-stage variables and a second-stage with continuous variables and a single binary variable per scenario (i.e., the export ban activation).

\begin{table}[!htb]  
\scriptsize
\vspace{-10pt}
\renewcommand{\arraystretch}{0.9}
\setlength{\tabcolsep}{2pt}
  \caption{Optimization Model Notation}
\begin{tabular}{ll}
\hline
\multicolumn{2}{l}{\textbf{Sets}}  \\
$I$ & supplier countries.\\ 
$J$ & potential countries to locate a manufacturing plant.\\
$\bar{J}$ & countries that ban exports, where $\bar{J} \subseteq J$.\\
$K$ & countries with drug demand, where $I\subseteq K$ and $J\subseteq K$.\\
$L$ & levels of manufacturing plant quality.\\
$\Omega$ & scenarios.\\
\multicolumn{2}{l}{\textbf{Parameters}} \\ 
$\varrho^{\omega}$ & probability of scenario $\omega \in \Omega$.\\
$c_{i}^{s}$ & unit cost of raw material purchased from supplier $i \in I$. \\ 
$c_{j}^{p}$ & unit drug production cost at plant $j \in J$. \\
$c_{jl}^{f}$& annualized fixed cost to produce drug at plant $j \in J$ with quality $l \in L$.\\
$c_{k{k^\prime}}^{t}$ & unit transportation cost from country $k \in K$ to country $k^{\prime} \in K$.\\
$c_{j}^{h}$ & unit excess raw material holding cost in plant $j \in J$.\\
$p_k$ & drug's unit price in country $k \in K$.\\
$q_{i}^{s}$ & drug capacity of supplier $i \in I$. \\
$q_{j}^{p}$ & drug capacity of plant $j \in J$.\\
$\eta_k$ & exogenous drug exports from country $k \in K$.\\
$\bar{g}$ & nominal global drug production.\\
$r$ & sensitivity threshold for the nominal global drug availability.\\
\multicolumn{2}{l}{\textbf{Random Parameters} (All associated with a scenario index $\omega$)} \\
\multicolumn{2}{l}{\textit{Operational Uncertainty}}\\
$\gamma_i^{s,\omega}$& proportion of capacity of supplier $i \in I$ that is available (operational strain).\\
$\gamma_{jl}^{p,\omega}$& proportion of capacity of plant $j \in J$ with quality $l \in L$ that is available (operational strain).\\
$d_k^{\omega}$& company-experienced drug demand in country $k \in K$.\\
\multicolumn{2}{l}{\textit{Disruption Risks}}\\
$\xi_k^{nd,\omega}$& binary parameter that equals 1 if capacity in country $k \in K$ is not disrupted by natural disasters.\\
$\xi_{i}^{qs,\omega}$& binary parameter that equals 1 if supplier $i \in I$ does not experience a quality-type disruption.\\
$\xi_{jl}^{qp,\omega}$& binary parameter that equals 1 if plant $j \in J$ with quality $l \in L$ does not experience a quality-type disruption.\\
\multicolumn{2}{l}{\textbf{Decision Variables}} \\
\multicolumn{2}{l}{\textit{First-Stage}}\\
$y_{jl}$ & binary variable that equals 1 if plant $j \in J$ with quality $l \in L$ is selected and 0 otherwise.\\
\multicolumn{2}{l}{\textit{Second-Stage} (All associated with a scenario index $\omega$)}\\
$u_{ij}^{\omega}$ & raw material purchased from supplier $i \in I$ and sent to plant $j \in J$. \\
$x_{j}^{\omega}$ & planned drug production in plant $j \in J$. \\
$\tau_{j{j^\prime}}^{\omega}$ & raw materials transshipped from  plant $j \in J$ to plant $j^{\prime} \neq j \in J$. \\
$e_{j}^{\omega}$ & excess raw materials in plant $j \in J$. \\ 
$v_{jk}^{\omega}$ & drugs produced by plant $j \in J$ and sent to country $k \in K$. \\
$z^{\omega}$ & auxiliary binary variable that equals 0 if countries in set $\bar{J}$ impose export bans.\\
\hline
\end{tabular}\label{Table:P2Notation}
\vspace{-10pt}
\end{table}

Next, we discuss the disruption and operational uncertainty in more detail. There are two types of conventional capacity disruptions that represent long-term facility closures: quality disruptions and disruptions due to natural disasters (e.g., \cite{WSJ-Mosbergen-Quality2023, NBC-Lovelace-Tornado2023}). Both are modeled as a binary parameter where a value of 1 represents availability, and 0 represents a disruption. Quality disruptions may occur at suppliers ($\xi_{i}^{qs,\omega} = 0, i\in I$) or at plants with a selected quality level ($\xi_{jl}^{qp,\omega}=0, j\in J, l\in L$). The probability of a quality disruption at a plant $j\in J$ varies by the selected quality level $l\in L$ ($y_{jl}$) where higher quality plants experience fewer disruptions. Natural disaster-induced disruptions $(\xi_k^{nd,\omega})$ may occur at suppliers $k\in I\subseteq K$ and at plants $k\in J\subseteq K$ if they are in a country affected by a natural disaster. Natural disaster-induced disruptions are correlated, following the shared hazard and propagation approach: if a natural disaster occurs, it may propagate to countries with shared hazards, and facilities within affected countries may be disrupted. Conventional operational uncertainty is represented by operational strain at suppliers and plants (i.e., the proportion of capacity that is available at supplier $i$ ($\gamma_i^{s,\omega}$) and plant $j$ ($\gamma_{jl}^{p,\omega}$) with quality level $l$, for all $i\in I, j\in J, l\in L$) as well as company-experienced demand ($d_k^\omega$) in country $k\in K$. Further details are in the supplementary materials (S\ref{supp:stochastic}).

Export bans are then an induced, systemic response to realized conventional disruptions and operational uncertainty; the relationship is given in Inequality \eqref{eq:EB_x}.
\begin{equation} \label{eq:EB_x}
\sum_{j \in J}x_{j}^\omega + \sum_{k \in K}\xi_k^{nd,\omega}\eta_k < r\bar{g} \qquad \Rightarrow \qquad z^\omega = 0, \qquad \omega \in \Omega.   
\end{equation}
The left-hand side represents the projected worldwide drug supply in a scenario $\omega \in \Omega$. The first term is the company's planned production, $x_j^\omega$, at plant $j\in J$ based on the realized natural disasters, quality disruptions, and operational uncertainty. The second term represents total production of the drug from other companies manufacturing in country $k \in K$ that is not disrupted by a natural disaster; other companies' production $\eta_k$ is an exogenous parameter representing production that would nominally be exported from country $k\in K$ if exports are allowed from $k$. The right-hand side represents the threshold for nominal global drug availability below which would activate an export ban for countries in $\bar{J}$. The nominal global drug availability, $\bar{g}$, is the total drug production (company and others) under regular supply chain operations, i.e., when the model is optimized without disruptions (supplementary materials (S\ref{supp:nominal})). The threshold $r$ represents the sensitivity of countries in set $\bar{J}$ to activating export bans; a higher value indicates higher sensitivity. For example, a threshold of $r=0.8$ indicates that if projected drug supply is below 80\% of nominal supply, then there will be systemic bans.

We draw attention to two subtleties with respect to how demand and export bans are defined. Countries in set $\bar{J}$ are sensitive to the worldwide projected supply \eqref{eq:EB_x} that includes both the company's production and other companies' exogenous production; this reflects policymakers' focus on availability of the drug from \textit{any} manufacturer rather than from one in particular. In contrast, the company can sell up to a portion of the total demand, and this is expressed as the company-experienced demand, $d_k^\omega$, in country $k\in K$ in scenario $\omega\in \Omega$. We further note that export bans are activated based on projected vs. nominal production rather than projected vs. 0\% shortages because there may be production shortfalls even in the nominal case, due to the underlying economics of generic injectable drugs \citep{SABOGAL2025}.

\subsection{\emph{Model Formulation}} \label{sub:formulation-2stage}
Next, we present a two-stage stochastic program for the supply chain design of a generic injectable drug manufacturer considering conventional and systemic disruption risks \eqref{1-2stage}.

\allowdisplaybreaks
\begin{subequations}\label{1-2stage}
    \footnotesize
        \begin{align}
             \max \limits_{\mathbf{y}} & \left \lbrace - \sum_{j \in J} \sum_{l \in L} c_{jl}^{f}y_{jl} + \sum_{\omega \in \Omega}\varrho^{\omega}\mathbf{Q}^{\omega}\mathbf{(y)} \right \rbrace \label{OF1-2stage}\\
             \text{s.t.  } & \sum_{l \in L}y_{jl} \leq 1 \quad \forall j \in J \label{1.1-2stage}\\
             & y_{jl} \in \{0,1\} \quad \forall j \in J, l \in L \label{1.2-2stage}
        \end{align}
        \begin{align}
             \mathbf{Q}^{\omega}\mathbf{(y)}:= \max & \sum_{j \in J} \sum_{k \in K} (p_{k}-c_{j}^{p}-c_{jk}^{t})v_{jk}^{\omega} - \sum_{i \in I}\sum_{j \in J}(c_{i}^{s}+c_{ij}^{t})u_{ij}^\omega - \sum_{j \in J} \bigg (c_j^{h}e_j^{\omega} + \sum_{j^\prime \neq j \in J}c_{jj^\prime}^{t} \tau_{jj^\prime}^{\omega} \bigg ) \label{OF2-2stage}\\
             \text{s.t.  } & \sum_{j \in J}u_{ij}^\omega \leq q_{i}^{s}\gamma_i^{s,\omega}\xi_i^{nd,\omega}\xi_i^{qs,\omega} \quad \forall i \in I \label{2.1-2stage}\\
             & x_{j}^\omega \leq q_{j}^{p}\xi_j^{nd,\omega}\sum_{l \in L} \gamma_{jl}^{p,\omega}\xi_{jl}^{qp,\omega}y_{jl} \quad \forall j \in J \label{2.2-2stage}\\
             & x_{j}^\omega = \sum_{i \in I} u_{ij}^\omega \quad \forall j \in J \label{2.3-2stage}\\
             & \sum_{j \in J}x_{j}^\omega \leq \sum_{k \in K}d_{k}^\omega \label{2.4-2stage}\\
             & \sum_{k \in K}v_{jk}^\omega \leq q_{j}^{p}\xi_j^{nd,\omega}\sum_{l \in L} \gamma_{jl}^{p,\omega}\xi_{jl}^{qp,\omega}y_{jl} \quad \forall j \in J \label{2.5-2stage}\\
             & \sum_{k \in K}v_{jk}^{\omega} + e_j^{\omega} + \sum_{j^{\prime} \neq j \in J}\tau_{jj^{\prime}}^{\omega} = \sum_{i \in I}u_{ij}^{\omega}+\sum_{j^{\prime} \neq j \in J}\tau_{j^{\prime}j}^{\omega} \quad \forall j \in J \label{2.6-2stage}\\
             & \sum_{j^{\prime}\neq j \in J} \tau_{j^{\prime}j}^\omega \leq M_3 \sum_{l \in L}y_{jl} \quad \forall j \in J \label{2.12-2stage}\\
             & \sum_{j \in J}v_{jk}^\omega \leq d_k^\omega \quad \forall k \in K \label{2.7-2stage}\\
             & r\bar{g} - \left ( \sum_{j \in J} x_j^\omega + \sum_{k \in K}\xi_k^{nd,\omega}\eta_k \right ) \leq M_1(1-z^\omega) \label{2.8-2stage}\\
             & v_{jk}^\omega \leq M_2z^\omega \quad \forall j \in \bar{J}, k \neq j \in K\label{2.9-2stage}\\
            & u_{ij}^\omega, x_j^\omega, e_j^\omega, \tau_{jj^{\prime}}^{\omega}, v_{jk}^\omega \geq 0 \quad \forall i \in I, j \in J, j^{\prime} \neq j \in J, k \in K \label{2.10-2stage}\\
            & z^\omega \in \{0,1\}  \label{2.11-2stage}
        \end{align} 
\end{subequations}
The objective of the model is to maximize the expected profit~\eqref{OF1-2stage}. This includes the fixed costs of plants and the expected profit across scenarios $\omega \in \Omega$ given by function $\mathbf{Q}^{\omega}\mathbf{(y)}$. The second-stage objective function \eqref{OF2-2stage} considers the costs to produce and ship the finished form of the drug, the costs to purchase raw materials and ship them to plants, and excess costs of raw material and transshipment, respectively. In the first stage, constraints \eqref{1.1-2stage} enforce that the drug can be produced in at most one plant in a given country, and \eqref{1.2-2stage} enforce the variable domains. 

In the second stage, raw materials purchases \eqref{2.1-2stage} and drug production plans \eqref{2.2-2stage} are limited by the available capacity, which are affected by conventional disruptions (natural disasters and quality disruptions) and operational strain. Raw materials are purchased and shipped to plants to meet production plans \eqref{2.3-2stage}. Production plans cannot exceed company-experienced demand \eqref{2.4-2stage}. Constraints \eqref{2.5-2stage} limit the actual drug production in each plant to the  capacity available (affected by conventional disruptions and operational strain). Drugs are produced using raw materials at plant $j \in J$ or obtained through transshipment between plants \eqref{2.6-2stage}. The transshipment of raw materials is allowed only between plants that are open \eqref{2.12-2stage}. Constraints \eqref{2.7-2stage} ensure that the company does not produce and sell more drugs than company-experienced demand. Export bans are expressed with constraints \eqref{2.8-2stage}--\eqref{2.9-2stage} representing the relationship given in \eqref{eq:EB_x}. That is, when the projected worldwide drug supply falls short of the threshold for the nominal global drug availability, the auxiliary binary variable $z^\omega$ is forced to be 0 \eqref{2.8-2stage}, and thereby, plants in countries in set $\bar{J}$ are prevented from exporting drugs. Otherwise, $z^\omega$ may be 0 or 1 without influence (see supplementary materials (S\ref{supp:sec:EBcon_detail}) for further details). Finally, constraints \eqref{2.10-2stage} and \eqref{2.11-2stage} enforce the variable domains. Note that we can set the big-M parameters as follows: $M_1=r\bar{g}$, $M_2=\max\{q_j^p \vert j \in J\}$, and $M_3=\sum_{i \in I}q_{i}^{s}$.

\subsection{\emph{Model Features and Assumptions}} \label{sub:assumptions}

First, the model only considers direct, first-order effects of export bans (i.e., the prohibition of international flows), as the focus of this paper is the study of systemic disruptions induced by conventional disruptions. Additional factors such as price increases (studied elsewhere, \cite{SABOGAL2025}) and second-order stochastic and competitive effects of other companies' response (unstudied elsewhere) are not included, e.g., other companies' production is modeled as exogenous. Focusing on the first-order effects for a particular company allows the paper to study the fundamental relationships of induced disruptions as well as explore decomposition approaches for the unique model structure. The model applies a common global threshold for triggering such bans ($r$); this could be relaxed by stratifying parameter $r$ by $j\in \bar{J}$. We apply export bans exclusively to flows of final products (rather than raw materials) because these are the most affected by export prohibitions \citep{WTO2022a,WTO2022b}. For simplicity, we assume that there is at most one supplier and one candidate plant location in each country, consistent with previous work \citep{SABOGAL2025}; this assumption could be relaxed by adding sets of suppliers and candidate plant locations in each country. Finally, conventional and induced disruptions represent long-term capacity interruptions. This follows the real-world rigidity of recovering pharmaceutical supply chains, and the probability distributions used to parameterize them represent these low-likelihood but high-impact events (see supplementary materials (S\ref{supp:Data})).

The model has fixed recourse and holds the complete recourse property because we allow shortages in the second stage. The recourse function $\mathbf{Q^\omega}(\mathbf{y}) < \infty$ for all first-stage decisions $\mathbf{y}$ and scenarios $\omega \in \Omega$ because of demand and capacity constraints. In addition, $\mathbf{Q^\omega}(\mathbf{y}) \geq 0$ since this is a maximization problem and there is always the option of not producing, i.e., $x_{j}^\omega=0$ for all $j \in J$, $u_{ij}^{\omega}=0$ for all $i \in I, j \in J$ and $v_{jk}^\omega=0$ for all $j \in J, k \in K$.
	
\section{Solution Methods} \label{s:methods}
The proposed model \eqref{1-2stage} is a two-stage stochastic mixed-binary program with binary decision variables in the first and second stages. Solving its deterministic equivalent formulation (also known as the extensive form) using advanced commercial solvers is computationally intractable when the number of scenarios $|\Omega|$ is large. To address this limitation, decomposition approaches can be used that iteratively solve smaller problems until convergence is reached. We use the unique structure of a single binary variable per scenario to
adopt three Benders-style decomposition approaches.

The core of the approach follows Benders decomposition. Program \eqref{1-2stage} is decomposed into a master problem, $(MP)$, that contains first-stage variables and constraints (model \eqref{Master_Singlecut}) and Benders subproblems, each with only second-stage variables and constraints for a given scenario ($\mathbf{Q}^{\omega}\mathbf{(y)}$). The variable $\theta^\omega$ is an upper-estimator of $\mathbf{Q}^{\omega}\mathbf{(y)}$ and is used to approximate $\mathbf{Q}^{\omega}\mathbf{(y)}$ through cuts. These cuts are stored in the set $\mathcal{O}$. We do not generate feasibility cuts because the problem has complete recourse (Section \ref{sub:assumptions}). We use the notation $\boldsymbol{\theta}$ to represent the vector of all $\theta^{\omega}$ for $\omega \in \Omega$.
\begin{subequations}\label{Master_Singlecut}
    \small
        \begin{align}
            (MP): \max \limits_{y} &  - \sum_{j \in J} \sum_{l \in L} c_{jl}^{f}y_{jl} + \sum_{\omega \in \Omega}\varrho^{\omega}\theta^{\omega} \\
             \text{s.t.  } & \eqref{1.1-2stage}-\eqref{1.2-2stage}\\
             &(\mathbf{y},\boldsymbol{\theta}) \in \mathcal{O} \label{generalcutMP}
        \end{align}
\end{subequations}
We follow the alternating strategy proposed by \cite{Angulo2016}. We first attempt to cut off an incumbent master solution using Benders (L-shaped) optimality cuts before trying to add other, more complex alternative cuts, such as integer L-shaped cuts, disjunctive cuts, or bilinear cuts.
Note that given an incumbent master solution ($\hat{\mathbf{y}},\hat{\theta}^{\omega}$), the linear relaxation of subproblem $\mathbf{Q}^{\omega}\mathbf{(\hat{y})}$ (denoted $\mathbf{Q_{LP}^{\omega}}\mathbf{(\hat{y})}$) provides an upper bound, i.e., $\mathbf{Q_{LP}^{\omega}}\mathbf{(\hat{y})} \geq \mathbf{Q}^{\omega}\mathbf{(\hat{y})}$. Then, a sufficient condition to reject ($\hat{\mathbf{y}},\hat{\theta}^{\omega}$) is $\hat{\theta}^{\omega} > \mathbf{Q_{LP}^{\omega}}\mathbf{(\hat{y})}$. If this condition is met, we add a Benders optimality cut of the form \eqref{benderscut} to set $\mathcal{O}$ of the master problem $(MP)$ to cut off this solution, $\mathbf{\hat{y}}$. Otherwise, we seek to generate an alternative cut. 

We next describe the Benders optimality cut that will be used throughout the methods. Subproblem $\mathbf{Q^{\omega}}\mathbf{(\hat{y})}$ is a mixed-binary program that contains a single binary variable ($z^\omega$). Relaxing the binary variable allows us to use duality to construct Benders optimality cuts \eqref{benderscut} to obtain 
$\mathbf{Q_{LP}^{\omega}}\mathbf{(\hat{y})}$, which is an upper bound for $\mathbf{Q^{\omega}}\mathbf{(\hat{y})}$. 
\begin{equation}\label{benderscut}
    \theta^{\omega} \leq \sum_{j \in J}\sum_{l \in L}\alpha_{jl}^{\omega}y_{jl}+\beta^{\omega}
\end{equation}
The coefficients $\beta^\omega$ and $\alpha_{jl}^{\omega}$ are calculated using the optimal dual variables $\varphi^{(\cdot),\omega}$ of $\mathbf{Q_{LP}^{\omega}}\mathbf{(\hat{y})}$ as follows, where $(\cdot)$ refers to the corresponding constraint in (\ref{1-2stage}).
\begin{small}
    \begin{align*}
       \beta^{\omega}  &=  \sum_{i \in I}\varphi_{i}^{\eqref{2.1-2stage},\omega} q_{i}^{s}\gamma_i^{s,\omega}\xi_i^{nd,\omega}\xi_i^{sq,\omega}+ \sum_{k \in K} \bigg ( (\varphi^{\eqref{2.4-2stage},\omega} + \varphi_k^{\eqref{2.7-2stage},\omega})d_k^\omega + \varphi^{\eqref{2.8-2stage},\omega}\eta_k\xi_k^{nd,\omega} \bigg ) 
    \end{align*}
    \vspace{-20pt}
    \begin{align*}
        \alpha_{jl}^{\omega} & = q_{j}^{p}\xi_j^{nd,\omega}\gamma_{jl}^{p,\omega}\xi_{jl}^{pq,\omega} \big (\varphi_{j}^{\eqref{2.2-2stage},\omega} + \varphi_{j}^{\eqref{2.5-2stage},\omega} \big) + \varphi_{j}^{\eqref{2.12-2stage},\omega}M_3
    \end{align*}
\end{small}
\subsection{\emph{Integer L-shaped Cut}} \label{sub:Int-Lshaped}
The integer L-shaped cut introduced by \cite{Laporte1993} is designed to solve two-stage stochastic programs with binary first-stage decision variables \citep{Angulo2016}. This cut is given by \eqref{intLshapedcut} and requires the optimal objective value of the mixed-integer subproblem $\mathbf{Q^{\omega}}\mathbf{(\hat{y})}$, as well as an upper bound $U^{\omega} \in \mathbb{R}$ on $\mathbf{Q^{\omega}(y)}$. 
\begin{equation}\label{intLshapedcut}
    \theta^{\omega} \leq \mathbf{Q^{\omega}(\hat{y})}+ \left (U^{\omega} - \mathbf{Q^{\omega}(\hat{y})} \right ) \left ( \sum_{(j,l):\hat{y}_{jl}=0} y_{jl}+ \sum_{(j,l):\hat{y}_{jl}=1} (1- y_{jl}) \right )
\end{equation}
\begin{equation}\label{eq:U}
    U^{\omega} = \sum_{k \in K}d_{k}^{\omega}\big (p_{k}-\min_{j\in J}\{c_{j}^{p}\}-\min_{i\in I}\{c_{i}^{s}\} \big )
\end{equation}

The upper bound $U^\omega$ is defined as the total profit obtained by satisfying all demand, having the lowest possible raw material and production costs, and incurring no logistical costs \eqref{eq:U}. Note this bound is valid for instances where $p_k \geq \min_{j\in J}\{c_{j}^{p}\}+\min_{i\in I}\{c_{i}^{s}\}, \forall k\in K$, including those we study in the computational experiments.

\begin{algorithm}
\caption{Alternating Integer L-shaped Cut Method}\label{Algorithm:AI}
\begin{algorithmic}[1]
\spacingset{1}
\small
\STATE \textbf{Input:} Incumbent solution $(\boldsymbol{\hat{\theta}}, \mathbf{\hat{y}})$ of $(MP)$
\FOR{$\omega \in \Omega$}
    \STATE Solve $\mathbf{Q_{LP}^{\omega}(\hat{y})}$ 
    \IF{$\hat{\theta}^{\omega} > \mathbf{Q_{LP}^{\omega}(\hat{y})}$}
        \STATE Add a Benders optimality cut \eqref{benderscut} to set $\mathcal{O}$ 
        \IF{optimal $z^{\omega} \in \{0,1\}$}
            \STATE Add an integer L-shaped cut \eqref{intLshapedcut} to set $\mathcal{O}$
        \ENDIF 
    \ELSE
        \STATE Solve $\mathbf{Q^{\omega}(\hat{y})}$
        \IF{$\hat{\theta}^{\omega} > \mathbf{Q^{\omega}(\hat{y})}$}
            \STATE Add an integer L-shaped cut \eqref{intLshapedcut} to set $\mathcal{O}$
        \ENDIF   
    \ENDIF
\ENDFOR
\end{algorithmic}
\end{algorithm}

Given an incumbent solution to the master problem ($\hat{\mathbf{y}},\hat{\theta}^{\omega}$), if $\hat{\theta}^{\omega} > \mathbf{Q^{\omega}(\hat{y})}$, then an integer L-shaped cut \eqref{intLshapedcut} can be added to set $\mathcal{O}$ of $(MP)$ to cut off that solution. Algorithm \ref{Algorithm:AI} presents a cutting plane method that alternates between Benders optimality cuts \eqref{benderscut} and integer L-shaped cuts \eqref{intLshapedcut} to solve the proposed model \eqref{1-2stage}. Note that if the auxiliary variable $z^\omega$ is binary in the optimal solution of $\mathbf{Q_{LP}^{\omega}}\mathbf{(\hat{y})}$, then an integer L-shaped cut \eqref{intLshapedcut} can be generated as a byproduct when attempting to add a Benders optimality cut. Numerical experiments indicate that, in this situation, it is more effective to add both the Benders and the integer L-shaped cuts rather than only one (see supplementary materials (S\ref{supp:comp_results})).

\subsection{\emph{Disjunctive Cut}} \label{sub:disjuctive}

The unique characteristic of the second-stage problem, i.e., a single binary variable per scenario $\omega$, motivates the study of an enumeration-style approach that may be inefficient for models that have more binary variables in the second-stage. We generate a disjunctive cut by first solving two linear subproblems per scenario $\omega$: (i) $\mathbf{Q_{LP}^\omega}\mathbf{(\hat{y})}$ with constraint $z^\omega \leq 0$ that enforces export bans, and (ii) $\mathbf{Q_{LP}^\omega}\mathbf{(\hat{y})}$ with constraint $z^\omega \geq 1$. We refer to (i) and (ii) as $\mathbf{Q_{LP}^\omega}(\mathbf{\hat{y}},z^\omega=0)$ and $\mathbf{Q_{LP}^\omega}(\mathbf{\hat{y}},z^\omega=1)$, respectively. Then, we use duality
to build a Benders optimality cut for each problem. Applying the disjunctive cut principle to these two cuts, we 
obtain a valid inequality (disjunctive cut) that can be passed to ($MP$) since it approximates the subproblem value function, $\mathbf{Q^\omega(y)}$. We follow the procedure presented by \cite{Sen-Sherali2006} to generate the disjunctive cuts.

The following presents the detailed walkthrough of the approach. The optimal solution for $\mathbf{Q^\omega}(\mathbf{\hat{y}})$ must be associated with at least one of the two nodes, $\mathbf{Q_{LP}^\omega}(\mathbf{\hat{y}},z^\omega=0)$ or $\mathbf{Q_{LP}^\omega}(\mathbf{\hat{y}},z^\omega=1)$. Then, disjunction \eqref{eq:disj} must hold:
\begin{equation}\label{eq:disj}
     \theta^{\omega} \leq \sum_{j \in J} \sum_{l \in L} \alpha_{0jl}^{\omega}y_{jl}+\beta_{0}^{\omega} \quad \vee \quad \theta^{\omega} \leq \sum_{j \in J} \sum_{l \in L} \alpha_{1jl}^{\omega}y_{jl}+\beta_{1}^{\omega},
\end{equation}
where coefficients $(\alpha_{0jl}^{\omega},\beta_0^{\omega})$ and $(\alpha_{1jl}^{\omega},\beta_1^{\omega}$) are calculated using the optimal dual variables of $\mathbf{Q_{LP}^\omega}(\mathbf{\hat{y}},z^\omega=0)$ and $\mathbf{Q_{LP}^\omega}(\mathbf{\hat{y}},z^\omega=1)$, respectively, as presented in \eqref{eq:B0}, \eqref{eq:B1} and \eqref{eq:alphaz}.
\begin{small}
    \begin{align}\label{eq:B0}
       \beta_0^{\omega} = \sum_{i \in I}\varphi_{i}^{\eqref{2.1-2stage},\omega} q_{i}^{s}\gamma_i^{s,\omega}\xi_i^{nd,\omega}\xi_i^{sq,\omega}+ \sum_{k \in K} \bigg ( (\varphi^{\eqref{2.4-2stage},\omega} + \varphi_k^{\eqref{2.7-2stage},\omega})d_k^\omega + \varphi^{\eqref{2.8-2stage},\omega}\eta_k\xi_k^{nd,\omega} \bigg )
    \end{align}
    \vspace{-20pt}
    \begin{align}\label{eq:B1}
       \beta_1^{\omega} = \sum_{i \in I}\varphi_{i}^{\eqref{2.1-2stage},\omega} q_{i}^{s}\gamma_i^{s,\omega}\xi_i^{nd,\omega}\xi_i^{sq,\omega}+ \sum_{k \in K} \bigg ( (\varphi^{\eqref{2.4-2stage},\omega} + \varphi_k^{\eqref{2.7-2stage},\omega})d_k^\omega + \varphi^{\eqref{2.8-2stage},\omega}\eta_k\xi_k^{nd,\omega} \bigg ) - \varphi^{(-z \leq -1),\omega}
    \end{align}
    \vspace{-20pt}
    \begin{align}\label{eq:alphaz}
        \alpha_{zjl}^{\omega} & = q_{j}^{p}\xi_j^{nd,\omega}\gamma_{jl}^{p,\omega}\xi_{jl}^{pq,\omega} \big (\varphi_{j}^{\eqref{2.2-2stage},\omega} + \varphi_{j}^{\eqref{2.5-2stage},\omega} \big) + \varphi_{j}^{\eqref{2.12-2stage},\omega}M_3 \quad \textit{for } z \in \{0,1\}
    \end{align}
\end{small}

Let the disjunctive set $\Pi(\omega)$ be the set of points $(\theta^{\omega},\mathbf{y})$ that are in at least one of the next two polyhedra.
\begin{small}
\[ \Pi(\omega)
=\left\{ (\theta^\omega, \mathbf{y}) \in \left\{ 
\begin{aligned}
    \theta^\omega & - \sum_{j \in J}\sum_{l \in L} \alpha_{0{jl}}^{\omega}y_{jl} \leq \beta_0^{\omega}  \\
    \sum_{l \in L} y_{jl} &\leq 1 \quad \forall j \in J \\
    y_{jl} &\leq 1 \quad \forall j \in J, l \in L \\
    y_{jl} &\geq 0 \quad \forall j \in J, l \in L
\end{aligned}
\right\} \bigcup \left\{
\begin{aligned}
    \theta^\omega & - \sum_{j \in J}\sum_{l \in L}  \alpha_{1{jl}}^{\omega}y_{jl} \leq \beta_1^{\omega}  \\
    \sum_{l \in L} y_{jl} &\leq 1 \quad \forall j \in J \\
    y_{jl} &\leq 1 \quad \forall j \in J, l \in L \\
    y_{jl} &\geq 0 \quad \forall j \in J, l \in L
\end{aligned}
\right\} \right\} \]
\end{small}

Applying the disjunctive cut principle \citep{DisjBook2018,Sen-Higle2005}, we can find non-negative multipliers $\lambda := [\lambda_{0}^{c}, \lambda_{1}^{c}, \lambda_{0j}^{m}, \lambda_{1j}^{m}, \lambda_{0jl}^{u}, \lambda_{1jl}^{u}]^T$ to form a valid inequality (disjunctive cut) for the convex hull of the union of the two polyhedrons. This cut is defined as \eqref{eq:disjuctive-cut-w} and provides an upper bound on $\mathbf{Q}^{\omega}\mathbf{(y)}$.
\begin{equation}\label{eq:disjuctive-cut-w}
    \bar{\delta}^{\omega}\theta^{\omega} + \sum_{j \in J} \sum_{l \in L} \bar{\alpha}_{jl}^{\omega} y_{jl} \leq \bar{\beta}^{\omega} 
\end{equation}
The cut coefficients $\bar{\delta}^{\omega}$, $\bar{\alpha}_{jl}^{\omega}$ and $\bar{\beta}^{\omega}$ are given by:
\begin{small} 
\begin{gather*}
    \bar{\delta}^{\omega}=\min\{\lambda_{0}^{c} \hspace{2pt}, \hspace{2pt} \lambda_{1}^{c}\} \\
    \bar{\alpha}_{jl}^{\omega}=\min \Big \{-\alpha_{0jl}^{\omega}\lambda_{0}^c  + \lambda_{0j}^{m} + \lambda_{0jl}^{u}  \hspace{2pt}, \hspace{2pt} -\alpha_{1jl}^{\omega}\lambda_{1}^c  + \lambda_{1j}^{m} + \lambda_{1jl}^{u}, \Big \} \quad \forall j \in J, l \in L \\
    \bar{\beta}^{\omega}=\max \Bigg \{\beta_{0}^{\omega}\lambda_{0}^c + \sum_{j \in J} \lambda_{0j}^{m} + \sum_{j \in J}\sum_{l \in L} \lambda_{0jl}^u \hspace{2pt}, \hspace{2pt} \beta_{1}^{\omega}\lambda_{1}^c + \sum_{j \in J} \lambda_{1j}^{m} + \sum_{j \in J}\sum_{l \in L} \lambda_{1jl}^u \Bigg \}
\end{gather*}
\end{small}

Algorithm \ref{Algorithm:D} in the supplementary materials (S\ref{supp:Algorithms}) presents a cutting plane method that alternates between Benders optimality cuts and disjunctive cuts to solve the proposed model \eqref{1-2stage}.
Note that the problem $\mathbf{Q_{LP}^\omega}(\mathbf{\hat{y}}, z^\omega=0)$ is always feasible, since constraint \eqref{2.8-2stage} is feasible for $x_j^{\omega} \geq 0$ for all $j \in J$ and \eqref{2.9-2stage} is feasible by setting $v_{jk}^\omega=0$ for $j \in \bar{J}, k \neq j \in K$. However, $\mathbf{Q_{LP}^\omega}(\mathbf{\hat{y}}, z^\omega=1)$ may be infeasible due to constraint \eqref{2.8-2stage}. The algorithm describes procedures for both cases.

\subsection{\emph{Bilinear Cut}} \label{sub:BilinearCuts}
Finally, we observe that the single binary variable, $z^\omega$, in subproblem $\mathbf{Q^{\omega}(\hat{y})}$ serves a similar role as the binary variable commonly included in the bilinear formulations of chance-constrained problems (e.g., \cite{Haghighat2018}). Specifically, when export bans are activated ($z^{\omega}=0$), the constraints that limit exports \eqref{2.9-2stage} are enforced, and when $z^\omega=1$, such constraints become redundant. In addition, the recourse function can be approximated using $z^\omega$ as a ``weight" that defines subproblem $\mathbf{Q^{\omega}(\hat{y})}$ as the convex combination of two linear subproblems $\mathbf{Q_{LP}^\omega}(\mathbf{\hat{y}}, z^\omega=0)$ and $\mathbf{Q_{LP}^\omega}(\mathbf{\hat{y}}, z^\omega=1)$. From this, duality can be used to build Benders optimality cuts, as presented by \eqref{eq:bil} and simplified to \eqref{eq:bil_2}.
\begin{equation}\label{eq:bil}
     \theta^\omega \leq \bigg ( \sum_{j \in J} \sum_{l \in L} \alpha_{0jl}^{\omega}y_{jl}+\beta_{0}^{\omega} \bigg )(1-z^{\omega})+ \bigg (\sum_{j \in J} \sum_{l \in L} \alpha_{1jl}^{\omega}y_{jl}+\beta_{1}^{\omega} \bigg) z^{\omega}
\end{equation}

\begin{equation}\label{eq:bil_2}
     \theta^\omega \leq \sum_{j \in J} \sum_{l \in L} \big (y_{jl}z^{\omega}(\alpha_{1jl}^{\omega}-\alpha_{0jl}^{\omega})+\alpha_{0jl}^{\omega}y_{jl} \big ) +(\beta_{1}^{\omega} - \beta_{0}^{\omega})z^{\omega}+\beta_{0}^{\omega}
\end{equation}

Note that the bilinear term $y_{jl}z^{\omega}$ corresponds to the product of two binary variables. This term is linearized by introducing the auxiliary variable $w_{jl}^{\omega}$ that captures the product of the binary variables ($y_{jl}z^{\omega}$) with the help of additional constraints in the master problem (see supplementary materials (S\ref{supp:Algorithms})). The reformulated linear cut is given by \eqref{master_bilinearcut} and is added iteratively in the algorithm. 

\begin{equation}\label{master_bilinearcut}
    \theta^\omega \leq \sum_{j \in J} \sum_{l \in L} \big ( w_{jl}^{\omega}(\alpha_{1jl}^{\omega}-\alpha_{0jl}^{\omega})+\alpha_{0jl}^{\omega}y_{jl} \big ) + (\beta_{1}^{\omega} - \beta_{0}^{\omega})z^{\omega} +\beta_{0}^{\omega}
\end{equation}

The master problem for solving model \eqref{1-2stage} using linearized bilinear cuts is given by $(MP^\prime)$, where $y_{jl}, \theta^\omega, z^\omega$ and $w_{jl}^{\omega}$ for $j \in J, l \in L, \omega \in \Omega$ are variables, and set $\mathcal{O^\prime}$ stores Benders optimality cuts and linearized bilinear cuts added iteratively. We use notation $\mathbf{z}$ and $\mathbf{w}$ to represent the vectors of all $z^{\omega}$ and $w_{jl}^{\omega}$ for $\omega \in \Omega, j \in J, l \in L$, respectively.
\begin{subequations}\label{Master2_Singlecut}
    \small
        \begin{align}
            (MP^\prime): \max \limits_{y} &  - \sum_{j \in J} \sum_{l \in L} c_{jl}^{f}y_{jl} + \sum_{\omega \in \Omega}\varrho^{\omega}\theta^{\omega} \\
             \text{s.t.  } & \eqref{1.1-2stage}-\eqref{1.2-2stage}\\
             & \eqref{lin_1}-\eqref{lin_4}\\
             &(\mathbf{y},\boldsymbol{\theta},\mathbf{z,w}) \in \mathcal{O^\prime} \label{generalcutMP2}\\
             & z^{\omega} \in \{0,1\} \quad \forall \omega \in \Omega
        \end{align}
\end{subequations}
Algorithm \ref{Algorithm:B} (supplementary materials (S\ref{supp:Algorithms})) presents a cutting plane method that alternates between Benders optimality cuts and linearized bilinear cuts to solve model \eqref{1-2stage}.

\section{Computational Experiments} \label{s:computational_results}
In this section, we report the computational results of applying the three cutting plane methods to solve model \eqref{1-2stage} and solving its extensive form with a commercial solver.

\subsection{\emph{Data Description and Test Instances}} \label{sub:data_exp}
The dataset used for the computational analyses and case study is based on vincristine (an oncology drug). It extends the dataset presented by \cite{SABOGAL2025} to include export prohibitions and natural disasters. It includes 11 suppliers of raw materials ($I$), 60 potential plant locations ($J$) each with three quality levels ($L$), and 179 countries with demand ($K$). The subset of countries that may implement export bans when nominal global drug availability falls below a threshold 
($\bar{J}$) consists of 18 countries that have implemented export prohibitions on pharmaceutical products \citep{WTO2022a,WTO2022b}. Each scenario is equally likely ($\varrho^\omega= 1/|\Omega|$ for all $\omega \in \Omega$). Complete details are in the supplementary materials (S\ref{supp:Data}).

We evaluate three sets of instances. Set I represents baseline case study parameters and is tested under export ban threshold values $r \in \{0.65, 0.8, 0.95\}$ and sample sizes $|\Omega| \in \{100, 500, 1000\}$, resulting in a total of nine instances. Sets II and III are designed to evaluate the computational effects of having export bans be activated more often and their implications (constraints \eqref{2.8-2stage} and \eqref{2.9-2stage}). In Set II, we double the number of countries that implement export bans ($|\bar{J}|=36$). In Set III, we increase the strength of the disruption propagation effect of natural disasters from a medium to a high level ($\zeta = 1300$ km), aligning with values reported in related studies (e.g., \cite{Li2018}). Since our goal is to test more challenging instances, we test Sets II and III with the export ban threshold $r = 0.95$ and sample sizes $|\Omega| \in \{100, 500, 1000\}$ representing a total of six additional instances. We generate five independent replications for each instance. 

\subsection{\emph{Implementation Details}} \label{sub:implementation_details}
All experiments are run on a single thread of a computational cluster (Palmetto) with eight cores, 50GB RAM, a time limit of 3600 seconds, and optimization solver Gurobi v11.0.0. For our cutting plane methods, we use the lazy constraint callback function available in Gurobi to add the cuts whenever a new master problem solution satisfying the integrality constraints is found. 

For all methods, a cut violation of $\epsilon = 10^{-5}$ is used. The initial master problem includes upper bounds for the $\theta^{\omega}$ variables, denoted by $U^{\omega}$ and computed as in \eqref{eq:U}. In Algorithm \ref{Algorithm:AI} (Alternating Integer L-shaped Cut), subproblem $\mathbf{Q^{\omega}(\hat{y})}$ is solved as a mixed-integer program ($z^\omega$ as a binary variable). In the other algorithms, Algorithm \ref{Algorithm:D} (Alternating Disjunctive Cut) and Algorithm \ref{Algorithm:B} (Alternating Bilinear Cut), two linear subproblems are solved per scenario: $\mathbf{Q_{LP}^\omega}(\mathbf{\hat{y}}, z^\omega=0)$ and $\mathbf{Q_{LP}^\omega}(\mathbf{\hat{y}}, z^\omega=1)$. 

In the proposed solution methods, the alternating strategy means first solving $\mathbf{Q_{LP}^{\omega}}\mathbf{(\hat{y})}$. If the variable $z^\omega$ is binary in the optimal solution of $\mathbf{Q_{LP}^{\omega}}\mathbf{(\hat{y})}$, Algorithm \ref{Algorithm:AI} adds two cuts to the master problem: a Benders optimality cut and an integer L-shaped cut. Algorithms \ref{Algorithm:D} and \ref{Algorithm:B} require one additional linear subproblem to be solved. For instance, if for a particular scenario $\hat{\omega}$, solving $\mathbf{Q_{LP}^{\hat{\omega}}}\mathbf{(\hat{y})}$ yields $z^{\hat{\omega}}=1$, and a Benders cut is not violated, then to generate a disjunctive or bilinear cut, we need to solve $\mathbf{Q_{LP}^{\hat{\omega}}}(\mathbf{\hat{y}}, z^{\hat{\omega}}=0)$. 

A computational issue can arise when generating the cuts in Algorithm \ref{Algorithm:D} if the resulting coefficient $\bar{\delta}^{\omega}=0$. In this case, a disjunctive cut is not generated for scenario $\omega$, and the algorithm proceeds to the next scenario. To avoid numerical issues, cut coefficients smaller than $10^{-6}$ are set to zero.

\subsection{\emph{Algorithms' Performance}} \label{sub:numerical_results}
In this section, we assess the computational performance of the three previously presented cutting plane methods (Algorithm~\ref{Algorithm:AI} to Algorithm~\ref{Algorithm:B}) and the performance of the extensive form solved by Gurobi directly using its default setting. 

Table \ref{table:computational_multicut} presents a summary of the computational results for baseline instances (Set I) across all methods. For each method, the table reports the average solution time (s) and the average number of branch-and-bound nodes explored of the instances solved to optimality within the time limit. These are under columns ``T" and ``Nds", respectively. The column ``Solved" indicates the number of instances (out of five) that were solved to optimality within the time limit. The column ``Gap" reports the average gap of the instances that were not solved to optimality.

\begin{table}[ht!]
\renewcommand{\arraystretch}{1}
\scriptsize
\setlength{\tabcolsep}{1.7pt} 
\caption{Computational performance for Set I (baseline case study)}
\begin{tabular}{c|c|cccc|cccc|cccc|cccc}
    \hline
    \multicolumn{2}{c|} {\textbf{Set I}} & \multicolumn{4}{c|}{\textbf{Alt. Integer L-shaped}} & \multicolumn{4}{c|}{\textbf{Alt. Disjunctive}} & \multicolumn{4}{c|}{\textbf{Alt. Bilinear}} & \multicolumn{4}{c}{\textbf{Extensive Form}} \\
    \hline
    $r$ & \textbf{$|\Omega|$} & \textbf{T} & \textbf{Nds} & \textbf{Solved} & \textbf{Gap} & \textbf{T} & \textbf{Nds} & \textbf{Solved} & \textbf{Gap} & \textbf{T} & \textbf{Nds} & \textbf{Solved} & \textbf{Gap}  & \textbf{T} & \textbf{Nds} & \textbf{Solved} & \textbf{Gap} \\
    \hline
        \multirow{3}{*}{0.65} & 100 & 39 & 163 & 5 & 0.0\% & 65 & 172 & 5 & 0.0\% & 245 & 612 & 5 & 0.0\% & 815 & 100 & 5 & 0.0\% \\
        & 500 & 360 & 430 & 5 & 0.0\% & 380 & 299 & 5 & 0.0\% & 928 & 666 & 5 & 0.0\% & - & - & 0 & - \\
        & 1000 & 675 & 366 & 5 & 0.0\% & 683 & 203 & 5 & 0.0\% & 1387 & 567 & 5 & 0.0\% & - & - & 0 & - \\
        \hdashline 
        \multirow{3}{*}{0.80} & 100 & 39 & 163 & 5 & 0.0\% & 67 & 172 & 5 & 0.0\% & 255 & 607 & 5 & 0.0\% & 754 & 104 & 5 & 0.0\% \\
        & 500 & 360 & 430 & 5 & 0.0\% & 371 & 300 & 5 & 0.0\% & 837 & 715 & 5 & 0.0\% & - & - & 0 & - \\
        & 1000 & 679 & 366 & 5 & 0.0\% & 789 & 221 & 5 & 0.0\% & 1493 & 568 & 5 & 0.0\% & - & - & 0 & - \\
        \hdashline 
        \multirow{3}{*}{0.95} & 100 & 108 & 389 & 5 & 0.0\% & 103 & 230 & 5 & 0.0\% & 302 & 633 & 5 & 0.0\% & 1365 & 126 & 5 & 0.0\% \\
        & 500 & 698 & 702 & 5 & 0.0\% & 932 & 524 & 5 & 0.0\% & 1116 & 630 & 5 & 0.0\% & - & - & 0 & - \\
        & 1000 & 1054 & 463 & 4 & 2.3\% & 1476 & 446 & 5 & 0.0\% & 1935 & 632 & 5 & 0.0\% & - & - & 0 & - \\
    \hline
\end{tabular}
\label{table:computational_multicut}
\end{table}

Our results indicate that the proposed decomposition methods are highly effective for solving problem \eqref{1-2stage}. In contrast, solving the extensive form with a commercial solver fails to yield optimal solutions for instances with a large number of scenarios (500 or 1000) within the one-hour time limit. In these instances, no incumbent solution is found. All 45 test instances are solved to optimality by the three proposed cutting plane methods, with the exception of one instance using the Alternating Integer L-shaped Cut Method, which ends with a gap of 2.3\%.

The Alternating Integer L-shaped Cut Method (Algorithm \ref{Algorithm:AI}) reports the best average computational time, followed by the Alternating Disjunctive Cut Method (Algorithm \ref{Algorithm:D}). Only in one instance ($r=0.95, |\Omega|=100$), Algorithm \ref{Algorithm:D} slightly outperforms Algorithm \ref{Algorithm:AI}, with an average solution time of 4.4 s faster. Additionally, for the instance ($r=0.95, |\Omega|=1000$), Algorithm \ref{Algorithm:D} presents the lowest computational times in two out of the five instances.

Algorithm \ref{Algorithm:AI} is, on average, 20.1 times faster (standard deviation 7.1) than solving the extensive form for the case of $|\Omega|=100$. It is surprising that Algorithm \ref{Algorithm:D} presents better average computational time than Algorithm \ref{Algorithm:B}, despite requiring the solution of up to three linear optimization problems for each scenario to generate a disjunctive cut. This is in contrast to Algorithm \ref{Algorithm:B}, which requires solving up to two linear optimization problems to generate a bilinear cut every time a feasible binary $\mathbf{y}$ solution is encountered. 
This result can be due to the larger size of the master problem in Algorithm \ref{Algorithm:B}, since the additional variables and constraints introduced by the linearization and the cut space $(z^{\omega}, y, w^{\omega}, \theta^{\omega})$. 

The computational performance for Sets II and III is presented in Table \ref{table:computational_ii_iii}. Both sets represent more challenging problem instances. Set II considers an increased number of countries that ban exports ($|\bar{J}|=36$). Set III considers natural disasters with a high strength of propagation effect ($\zeta = 1300$ km).

\begin{table}[ht!]
\renewcommand{\arraystretch}{1}
\scriptsize
\setlength{\tabcolsep}{1.5pt} 
\caption{Computational performance for Sets II ($|\bar{J}|=36$) and III (high propagation)}
\begin{tabular}{c|c|cccc|cccc|cccc|cccc}
    \hline
    \multicolumn{2}{c|} {$r=0.95$} & \multicolumn{4}{c|}{\textbf{Alt. Integer L-shaped}} & \multicolumn{4}{c|}{\textbf{Alt. Disjunctive}} & \multicolumn{4}{c|}{\textbf{Alt. Bilinear}} & \multicolumn{4}{c}{\textbf{Extensive Form}} \\
    \hline
    \textbf{Set} & \textbf{$|\Omega|$} & \textbf{T} & \textbf{Nds} & \textbf{Solved} & \textbf{Gap} & \textbf{T} & \textbf{Nds} & \textbf{Solved} & \textbf{Gap} & \textbf{T} & \textbf{Nds} & \textbf{Solved} & \textbf{Gap}  & \textbf{T} & \textbf{Nds} & \textbf{Solved} & \textbf{Gap} \\
    \hline
        \multirow{3}{*}{II} & 100 & 342 & 932 & 5 & 0.0\% & 692 & 629 & 5 & 0.0\% & 754 & 707 & 5 & 0.0\% & 1411 & 70 & 4 & 1.0\% \\
        & 500 & 3211 & 1500 & 2 & 3.1\% & - & - & 0 & 5.2\% & - & - & 0 & 1.2\% & - & - & 0 & - \\
        & 1000 & - & - & 0 & 6.8\% & - & - & 0 & 32.8\%* & - & - & 0 & 9.1\% & - & - & 0 & - \\
        \hdashline 
        \multirow{3}{*}{III} & 100 & 72 & 173 & 5 & 0.0\% & 125 & 398 & 5 & 0.0\% & 383 & 499 & 5 & 0.0\% & 1485 & 151 & 5 & 0.0\% \\
        & 500 & 743 & 506 & 5 & 0.0\% & 1384 & 621 & 5 & 0.0\% & 1099 & 607 & 5 & 0.0\% & - & - & 0 & - \\
        & 1000 & 1509 & 470 & 5 & 0.0\% & 1656 & 466 & 5 & 0.0\% & 2986 & 616 & 4 & 0.5\% & - & - & 0 & - \\
    \hline
\end{tabular}
\scriptsize *Average of three instances since the other two reported unbounded gaps
\label{table:computational_ii_iii}
\end{table}

Results in Table \ref{table:computational_ii_iii} demonstrate the increased complexity of decision-making under high levels of disruption risks. Among all implemented methods, Algorithm \ref{Algorithm:AI} (Alternating Integer L-shaped Cut) remains the most effective. It solves all 15 instances in Set III to optimality within one hour, as well as seven out of 15 instances in Set II. None of the methods are able to solve any instances from the most challenging instance to optimality within the time limit (Set II, $|\Omega| = 1000$, $r=0.95$); however, Algorithm \ref{Algorithm:AI} achieves the lowest average optimality gap in that case. In addition, for Set II with a high number of scenarios (500, 1000), Algorithm \ref{Algorithm:B} (Alternating Bilinear Cut) reports lower average gaps than Algorithm \ref{Algorithm:D} (Alternating Disjunctive Cut). For Set III, similar to Set I, Algorithm \ref{Algorithm:D} exhibits better computational performance than Algorithm \ref{Algorithm:B}. Considering the 64 instances that are solved to optimality by both Algorithms~\ref{Algorithm:AI} and \ref{Algorithm:D} across the three instance sets, Algorithm~\ref{Algorithm:AI} is, on average, 1.6 times faster, with a standard deviation of 1.1.

Next, we discuss additional implementations and combinations that are tested. Detailed results are presented in the supplementary materials (S\ref{supp:comp_results}). We evaluate the single-cut versions of Algorithms \ref{Algorithm:AI} to \ref{Algorithm:B} using Set I (baseline case study) under sensitivity threshold values $r \in \{0.65, 0.8, 0.95\}$ and sample sizes $|\Omega| \in \{100, 500, 1000\}$. In our results, among the single-cut methods, Algorithm \ref{Algorithm:AI} (Alternating Integer L-shaped Cut) achieves the best average computational time, followed by Algorithm \ref{Algorithm:D} (Alternating Disjunctive Cut). Additionally, the single-cut methods exhibit higher computational times compared to their multi-cut counterparts.

We also assess two additional cutting plane methods that combine three types of cuts. We consider the Alternating Integer L-shaped + Disjunctive Cuts and the Alternating Integer L-shaped + Bilinear Cuts. In these new methods, every time a disjunctive cut (bilinear cut) is added, an integer L-shaped cut is also added. This integration eliminates the need to solve the mixed-integer subproblem $\mathbf{Q^{\omega}(\hat{y})}$, as the procedures for generating disjunctive and bilinear cuts compute $\mathbf{Q^{\omega}(\hat{y})}$ through linear subproblems $\mathbf{Q_{LP}^\omega}(\mathbf{\hat{y}}, z^\omega=0)$ and $\mathbf{Q_{LP}^\omega}(\mathbf{\hat{y}}, z^\omega=1)$. We test these two additional algorithms with Set I under geopolitical threshold $r=0.95$ and sample sizes $|\Omega| \in \{100, 500, 1000\}$. Our results suggest that Algorithms \ref{Algorithm:AI} to \ref{Algorithm:B} are more effective than the combination of integer L-shaped cuts with disjunctive or bilinear cuts. This may be due to the larger size of the master problem, since for each scenario, two cuts are added instead of one, which also affects the scalability of the method as the number of scenarios grows. 

Finally, we implement feasibility cuts in Algorithm \ref{Algorithm:B} when $\mathbf{Q_{LP}^\omega}(\mathbf{\hat{y}}, z^\omega=1)$ is infeasible. This variation of Algorithm \ref{Algorithm:B} does not exhibit improved computational performance.

\section{Case Study} \label{s:value_risks}
Next, we evaluate the real-world implications of planning for systemic export bans. 

\subsection{\emph{Impact of Systemic Disruptions on Global Supply Chains}} \label{sub:SC_results}
For these analyses, we consider one sample of 1000 scenarios for each instance in Sets I, II, and III, all solved to optimality using the Alternating Integer L-Shaped Method (Algorithm \ref{Algorithm:AI}). For each instance set, we observe that export bans are not triggered when the sensitivity to projected global drug supply is $r=0.8$ or below ($r=0.65$), and the results are the same. When there is a higher sensitivity to projected supply shortages (here, $r=0.95$), there may be export bans. Accordingly, we compare the two cases $r=0.8$ and $r=0.95$, below.

\begin{figure}[ht!]
    \centering
    \includegraphics[width=15cm]{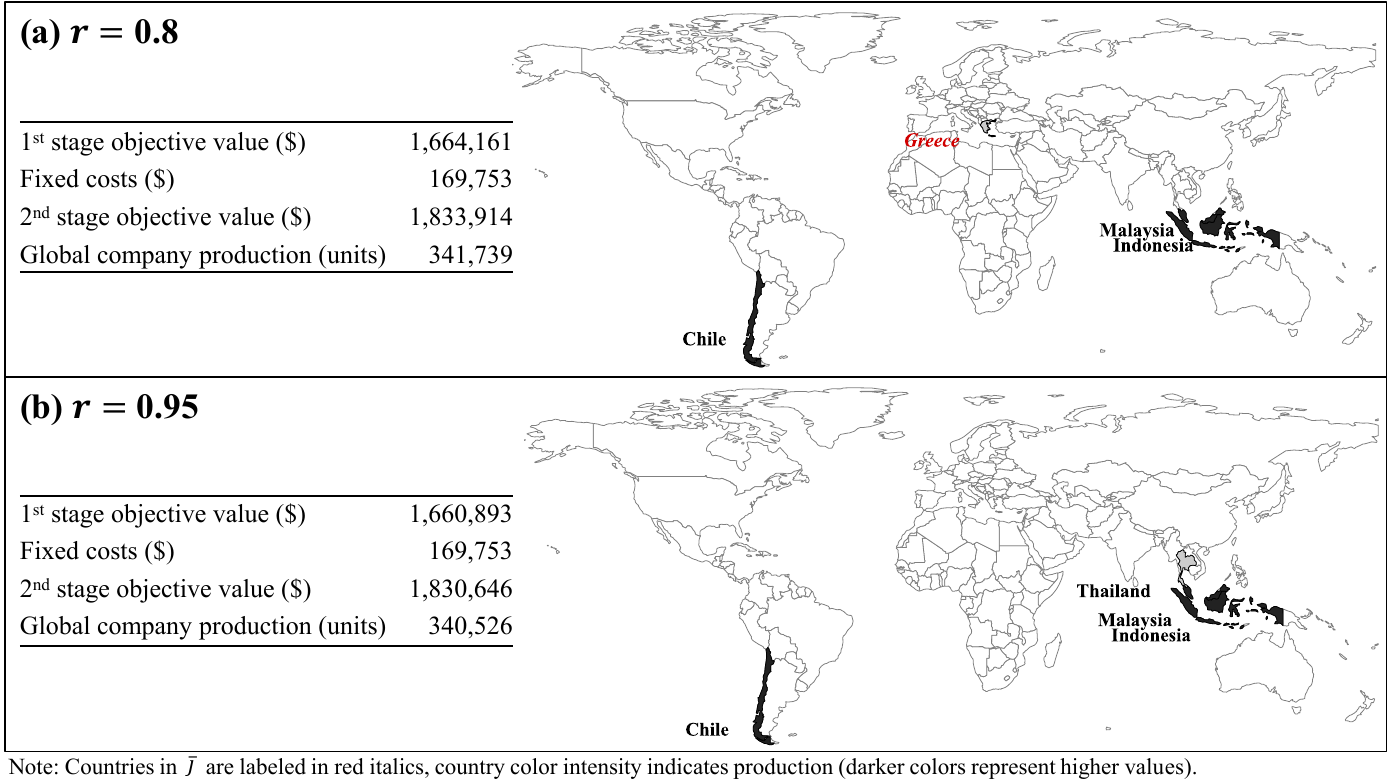}
    \caption{Plant locations for Set I instances with (a) $r=0.8$ and (b) $r=0.95$.}
    \label{fig:SetI.pdf}
    \vspace{-10pt}
\end{figure}

The results for the baseline instances (Set I) are shown in Figure \ref{fig:SetI.pdf}. Both supply chains consist of four plants at the lowest quality level, with three plants primarily serving demand, and one acting as a back-up. The primary plants are not at risk of imposing export bans (not in set $\bar{J}$) and are the same in both cases (Indonesia, Malaysia and Chile). These represent the three countries with the lowest production and fixed costs. There is a difference in the back-up plants between the two cases. Under moderate geopolitical risk ($r=0.8$), the back-up plant location, Greece is at risk of imposing export bans, whereas the back-up plant when there is higher sensitivity ($r=0.95$), Thailand is not at risk of imposing bans. Both Greece and Thailand have the same fixed costs, suggesting that planning for strain does not necessarily have higher strategic (first-stage) costs. Nonetheless, higher sensitivity to projected shortages do have minor negative effects on expected profits (0.2\% reduction), units sold (0.4\% decrease), and expected global shortages (0.3 absolute percentage points (pp) increase). The instances took 938 s ($r=0.8$) and 4,929 s ($r=0.95$) to solve.  

Next, we discuss the results for Set II  where more countries are at risk of imposing export bans ($|\bar{J}|=36$). A risk threshold of 0.95 vs. 0.8 reduces expected profits by 0.8\% (\$12,648) and increases fixed costs by 0.9\% (\$1,604), while the number of units sold remains unchanged. With $r=0.95$, all four selected countries (Greece, Malaysia, Indonesia and India) are in set $\bar{J}$. The backup plant is located in India, whose production costs are slightly higher than in Greece but still low (the seventh lowest production costs of all countries in $J$). 
To assess the challenge of a larger set of countries that may ban exports (Set II vs. Set I, $r=0.95$), we observe a 0.6\% (\$9,380) reduction in expected profits and a 0.9\% (\$1,604) increase in fixed costs. Surprisingly, the units sold increase by 0.4\% (\$1,309), and the expected global shortage decreases by 0.3 pp. Computationally, the model took 2.55 times longer to solve (17,506 s vs. 4,929 s). With $r=0.8$, it took 581 s to solve.

Finally, we discuss the results for Set III, which represents contexts with high propagation effects of natural disasters ($\zeta = 1300$ km) and a moderate number of countries imposing export bans ($|\bar{J}|=18$). Under both thresholds (0.8 and 0.95), the obtained supply chains are formed by Chile, Malaysia, Indonesia and Russia. None represent a threat of export bans (none belong to the set $\bar{J}$), and Russia has not experienced major natural disasters that completely disrupt manufacturing capabilities \citep{EM-DAT2024}. The performance of this supply chain does not substantially change compared to baseline disaster risks (Figure \ref{fig:SetI.pdf}). These instances took 849 ($r=0.8$) and 927 s ($r=0.95$) to solve.

\subsection{\emph{Value of Considering Disruption Risks}} \label{sub:VDR}
In this section, we quantify the value of integrating different types of disruption uncertainty into global supply chain design in two mis-specified settings: baseline that does not account for any disruptions and partial that accounts for some but not all disruptions. 

To capture the baseline value of planning for disruptions, we first solve the nominal model that does not include them (presented in supplementary materials (S\ref{supp:nominal})). We evaluate the resulting first-stage decisions, $\mathbf{\bar{y}}$, by fixing them as the first-stage decisions of \eqref{1-2stage} and solving with the uncertainty set that includes all disruptions, across different instances.

Table \ref{table:vdr} reports the benefits of planning for disruption risks vs. the mis-specified baseline. Across all considered instances, planning against disruptions leads to an increase in expected profits. For example, in the more challenging instance, Set II ($\bar{J}=36$) under high sensitivity threshold ($r=0.95$), the company achieves a $7.7\%$ increase in expected profits. Moreover, planning for disruptions also results in substantial reductions in the expected global drug shortage, explained by improved resilience. When planning for disruptions, four plants are selected, instead of three, enhancing the company's flexibility and ability to respond to disruptions.

\begin{table}[h]
\centering
\footnotesize
\caption{Value of considering disruption risks in supply chain design}
\begin{tabular}{ll|cc}
\hline
\multicolumn{2}{c|} {\textbf{Instance}} & \multicolumn{2}{c}{\textbf{Value of Planning for Disruption Risks}} \\
\textbf{Set} & $r$ & \textbf{Change in Expected Profits} & \textbf{Change in Expected Global Shortage} \\
\hline
I   & 0.95          & \$40,777 (2.5\%)    & -22.6\% (-6.0 pp)\\
    & 0.80 and 0.65 & \$17,002 (1.0\%)    & -19.2\% (-4.8 pp)\\
\hdashline
II  & 0.95          & \$118,327 (7.7\%)   & -30.4\% (-8.8 pp)\\
III & 0.95          & \$45,552 (2.8\%)    & -24.1\% (-6.4 pp)\\
\hline
\end{tabular}
\\ \tiny Note: pp, absolute percentage points
\label{table:vdr}
\end{table}

Next, we evaluate the value of planning for all three types of disruption risks (quality failures, $Q$, natural disasters, $ND$, and export bans, $EB$) against four partial settings that only consider a subset. The latter are called reference cases. The comparator with all disruption risks (model \eqref{1-2stage}) applies a geopolitical threshold of $r=0.95$ and uses Instance Set I. Case \textit{Q+EB} represents a model that only includes quality failures and export bans. This situation is implemented by fixing the natural disaster stochastic parameter $\xi_k^{nd,\omega}=1$ for all $k \in K, \omega \in \Omega$ and the sensitivity threshold to $r=0.95$. Case \textit{ND+EB} represents natural disasters and export bans only. This is done by setting the manufacturing quality failure stochastic parameter $\xi_{jl}^{qp,\omega}=1$ for all $j \in J, l \in L, \omega \in \Omega$ and the sensitivity threshold to $r=0.95$. Case \textit{Q} represents only quality failure disruptions. This is modeled by setting $\xi_k^{nd,\omega}=1$ for all $k \in K, \omega \in \Omega$ and a very low sensitivity threshold of $r=0.1$. The latter value means that export bans are practically negligible, and results confirm that there are no export bans induced in this case. Finally, Case \textit{ND} represents only natural disaster disruptions. This is represented by fixing $\xi_{jl}^{qp,\omega}=1$ for all $j \in J, l \in L, \omega \in \Omega$ and setting $r=0.1$, again rendering no export bans.

\begin{table}[h]
\centering
\footnotesize
\caption{Value of considering conventional and systemic disruption risks}
\begin{tabular}{ll|cc}
\hline
\textbf{Ref. Case} &\textbf{Comparator} & \textbf{Change in Expected Profits} & \textbf{Change in Expected Global Shortage} \\
\hline
 $Q+EB$ & All & \$0 (0.0\%)          & 0.0\% (0.0 pp)\\
 $ND+EB$ & All & +\$18,522 (+1.1\%)    & -18.2\% (-4.5 pp)\\
 $Q$     & All & +\$1,011 (+0.1\%)     & +0.3\% (+0.1 pp)\\
 $ND$    & All & +\$40,777 (+2.5\%)    & -22.6\% (-6.00 pp)\\
\hline
\end{tabular}
\\ \tiny Note: ref., reference; pp, absolute percentage points; +, increase; -, decrease; all, $Q+ND + EB$
\label{table:vmdr}
\end{table}

The value of considering all three disruptions is presented in Table \ref{table:vmdr}. Our results suggest that when the company plans for all three risks, it achieves the same expected profit and global shortage as when it plans for quality failures and export bans (Case \textit{Q+EB}) alone. This occurs because the plans with and without considering natural disasters result in identical supply chain configurations, with plants located in Indonesia, Malaysia, Thailand, and Chile. This finding suggests that, for the case study considered, natural disasters have limited influence on the company's strategic decisions. Under this four-plant structure, the sampled natural disaster events do not occur simultaneously across these countries, and a single backup facility is sufficient to mitigate capacity disruptions. In contrast, including quality disruptions (all vs. $ND+EB$ and all vs. $ND$), increases expected profits and decreases expected global shortages. Furthermore, planning for export bans (all vs. $Q$ and all vs. $ND$) also increases expected profits and decreases expected global shortages.

\section{Conclusion} \label{s:conclusion}

In this article, we introduce a global supply chain design model that integrates conventional and systemic disruption risks to produce a mixed-integer program with a pure binary first-stage problem and a nearly continuous second-stage problem. It is the first to consider induced export bans in a supply chain design context, and its special structure motivates a computational study on cutting plane methods. Experiments demonstrate that the Alternating Integer L-shaped Cut Method presents the best computational performance, followed by the Alternating Disjunctive Cut Method. Moreover, all three methods exhibit far superior performance compared to the extensive form. 

We apply the model to a case study of a life-saving oncology drug supply chain and evaluate how systemic disruptions may affect global strategy. We observe that both companies and patients benefit when the risks of systemic disruptions (export bans) are incorporated into supply chain design, with the most benefit for patients, i.e., lower shortages. We also observe that planning for quality disruptions and export bans capture most of the value of modeling disruptions, with the inclusion of correlated natural disasters having marginal effects on outcomes. This reflects our real-world observations as well, with quality issues more likely to cause drug shortages than natural disasters. It is interesting that despite this, all selected plants are the lowest quality level, suggesting that the marginal costs of increasing quality may not be worth the opportunity costs of quality-related disruption risks. We note however that if natural disaster frequency and severity increase, these disruptions may begin to have a more pronounced impact on supply chain design. Broadly, we observe that the effects of systemic disruptions are most influential when the countries' sensitivity to imposing them are high $r=0.95$; this suggests the utility of this modeling approach may increase in riskier environments. 

Future research in systemic, induced disruptions could focus on generating stronger disjunctive cuts in Algorithm \ref{Algorithm:D}. 
This algorithm generally exhibits the second-best computational time, with some cases achieving the best performance. Further, as systemic, global disruptions become more common, future work could also consider addressing key modeling limitations (see Section \ref{sub:assumptions}). In particular, a study of competitive dynamics between generic pharmaceutical companies under the risk of bans could be worthwhile.

\bibliographystyle{chicago}
\spacingset{1}
\bibliography{main-ref}

\newpage

\setcounter{section}{0}
\setcounter{figure}{0}
\setcounter{table}{0}
\setcounter{page}{1}

\input{Supp}

\bibliographystyleSupp{chicago}
\spacingset{1}
\bibliographySupp{supp-ref}

\end{document}

%% file: Supp.tex
\spacingset{1.5} 

\begin{center}
     \Large \bf {Supplementary Materials for ``When Shortages Lead to Export Restrictions: A Computational Study''}
     \end{center}

Further details on the capacity-related model parameters are presented in Section \ref{supp:stochastic}. The nominal model (that does not include disruptions) is in Section \ref{supp:nominal}. A walkthrough of the export ban constraints is presented in Section 
\ref{supp:sec:EBcon_detail}. Details on the data used in the case example are in Section \ref{supp:Data}. Additional algorithmic details are presented in Section \ref{supp:Algorithms}. Finally, additional results related to computational performance are in Section \ref{supp:comp_results}.

\section{Description of Capacity-Related Parameters} \label{supp:stochastic} 

In this section, we provide further details on the operational and disruption uncertainty. Recall that yield-style uncertainty is modeled as a proportion of total capacity available. For disruption uncertainty, a value of 1 indicates availability (non-disrupted), and a value of 0 indicates a disruption. 

\textit{Quality disruptions} $(\xi_{jl}^{qp,\omega},\xi_{i}^{qs,\omega})$ represent the complete unavailability of capacity due to quality failures in the plants and suppliers, respectively. These are associated with candidate plant locations $j\in J$, quality levels $l\in L$, and supplier countries $i\in I$. $\xi_{jl}^{qp,\omega} \sim {\rm Bernoulli}(\rho_{jl}^{qp})$ for all $j \in J$ and $l \in L$, where $\rho_{jl}^{qp}$ represents the probability that facility $j$ with quality $l$ is available. Similarly, $\xi_{i}^{qs,\omega} \sim {\rm Bernoulli}(\rho_{i}^{qs})$ for all $i \in I$, where $\rho_{i}^{qs}$ represents the probability that supplier $i$ is available.

\textit{Capacity strains} $(\gamma_{jl}^{p,\omega},\gamma_i^{s,\omega})$ are partial capacity interruptions of manufacturing plants and suppliers due to the regular uncertainty of the processes. Those are modeled with discrete random variables with probability mass functions $p(n)$, where $n$ represents capacity levels.

\textit{Natural disaster disruptions} $(\xi_k^{nd,\omega})$ represent the complete unavailability of the capacity of the company decision-maker and any exogenous producers in country $k\in k$. Specifically, $\xi_k^{nd,\omega}$ is a binary random parameter that equals 0 if production capacity in country $k \in K$ is disrupted by natural disasters that impacted country $k$. To generate (sample) capacity disruptions due to natural disaster disruptions in a particular country (i.e., the parameter used in the model), we step through three levels of their effects:
\begin{enumerate}
    \item Occurrence (epicenter) of natural disaster in country $k^\prime \in K$
    \item Propagation of natural disaster from epicenter country $k^\prime \in K$ to other country $k\in K$
    \item Risk of capacity disruption if country is affected by disaster (as an epicenter or via propagation)
\end{enumerate}

Level (1): the occurrence of natural disasters. Let $\mathcal{U}^{\omega}_{k^\prime}$ be a random parameter representing the non-occurrence of a natural disaster in country $k^\prime \in K$. $\mathcal{U}^{\omega}_{k^\prime} \sim {\rm Bernoulli}(1-\psi_{k^\prime})$, where $\psi_{k^\prime}$ is the probability that a natural disaster occurs in country $k^{\prime}$. If $\mathcal{U}^{\omega}_{k^\prime}$ equals 0, it means that country $k^\prime$ is the epicenter of a natural disaster. Consider $\bar{\mathcal{U}}$ as the set of countries that are epicenters of natural disasters, i.e., $\bar{\mathcal{U}} \coloneqq \{k^\prime \in K \vert \mathcal{U}^{\omega}_{k^\prime} = 0\}$.

Level (2): the propagation effect between an epicenter of a natural disaster, $k^\prime \in \bar{\mathcal{U}}$, and the countries with shared hazards with $k^\prime$. Let $H_{k^{\prime}}$ be the set of countries with common hazards with country $k^{\prime}$. This set is defined based on proximity in this paper. If a country $k \in K$ is not the epicenter of a natural disaster ($k \notin \bar{\mathcal{U}}$), and it does not have shared hazards with any country where a natural disaster occurred ($k \notin H_{k^\prime}$ for all $k^\prime \in \bar{\mathcal{U}}$), then capacity is not disrupted by a natural disaster in country $k\in K$, i.e., $\xi_{k}^{nd,\omega}=1$. 

Next, we consider the case when a country is not the epicenter of a natural disaster ($k \notin \bar{\mathcal{U}}$), but it has common hazards with another country where a natural disaster occurred ($k \in H_{k^\prime}$ for some $k^\prime \in  \bar{\mathcal{U}}$). Let $\mathcal{F}_{k{k^\prime}}^{\omega}$ be a binary random parameter that is equal to 1 if country $k \in H_{k^\prime}$ is not impacted by a natural disaster that occurred in country $k^\prime \in \bar{\mathcal{U}}$. $\mathcal{F}_{k{k^\prime}}^{\omega} \sim {\rm Bernoulli}(1-\rho_{k{k^\prime}}^{nd})$, where $\rho_{k{k^\prime}}^{nd}$ is the probability that country $k$ is impacted by a natural disaster that occurred in country $k^\prime$. This probability is estimated as $\rho_{k{k^\prime}}^{nd} = e^{\frac{-D_{k{k^\prime}}}{\zeta}}$, where $D_{k{k^\prime}}$ is the distance between countries $k$ and $k^\prime$, and $\zeta$ represents the strength of the disruption propagation effect. Consider $\bar{\mathcal{F}}$ as the set of countries that were impacted by natural disasters due to propagation effects, i.e, $k \in K\setminus \bar{\mathcal{U}}$ and $k \in \bigcup_{k^\prime \in \bar{\mathcal{U}}}H_{k^\prime}$ such that $\min\{\mathcal{F}_{k{k^\prime}}^{\omega} \vert k^\prime \in \bar{\mathcal{U}}$ and $k^{\prime} \in H_{k} \} = 0$. The approach to modeling the propagation effect follows standard methods used in other supply chain models (e.g., \cite{s-Lu-Ran2015} and \cite{s-Li2018}), where correlation is induced by shared hazard exposure and ripple effect-style disruptions. 

Level (3): risk of capacity disruption for countries impacted by natural disasters, i.e., $k \in \bar{\mathcal{U}}\cup\bar{\mathcal{F}}$. Plants in affected countries may (or may not) suffer disruptions in their production capacity. We sample whether production capacity is disrupted, using the conditional probability that manufacturing capacity in country $k$ is disrupted given that country $k$ was hit by natural disasters. 

\section{Nominal Model}\label{supp:nominal}

To calculate the nominal global drug availability, $\bar{g}$, we use a supply chain design model without disruptions (model \eqref{nominal-model}). Similar to the two-stage stochastic program \eqref{1-2stage}, $\bar{y}, \bar{u}, \bar{x}, \bar{\tau}, \bar{e}$, and $\bar{v}$ of model \eqref{nominal-model} are variables representing the location of manufacturing plants, purchases of raw materials, production plan, raw materials transshipment, excess of raw materials, and distribution of final products, respectively. $\varrho^{\bar{\omega}}$ corresponds to the probability of scenario $\bar{\omega} \in \bar{\Omega}$, given by the realization of operational uncertainties ($\gamma^{p,\bar{\omega}},\gamma^{s,\bar{\omega}},d^{\bar{\omega}}$). 

\begin{subequations}\label{nominal-model}
        \scriptsize
            \begin{align}
                 \max & - \sum_{j \in J} \sum_{l \in L} c_{jl}^{f}\bar{y}_{jl} + \sum_{\bar{\omega}\in \bar{\Omega}} \varrho^{\bar{\omega}} \left (\sum_{k \in K}\sum_{j \in J} (p_k- c_{j}^{p}-c_{jk}^{t}) \bar{v}_{jk}^{\bar{\omega}} - \sum_{i \in I}\sum_{j \in J}(c_{i}^{s}+c_{ij}^{t})\bar{u}_{ij}^{\bar{\omega}} - \sum_{j \in J} \bigg (c_j^{h}\bar{e}_j^{\bar{\omega}} + \sum_{j^\prime \neq j \in J}c_{jj^\prime}^{t} \bar{\tau}_{jj^\prime}^{\bar{\omega}} \bigg ) \right ) \label{OF4} \\                
                 \text{s.t.  } & \sum_{l \in L} \bar{y}_{jl} \leq 1 \quad \forall j \in J \label{4.1} \\
                 &\sum_{j \in J}\bar{u}_{ij}^{\bar{\omega}} \leq q_{i}^{s}\gamma_i^{s,\bar{\omega}} \quad \forall i \in I, \bar{\omega} \in \bar{\Omega} \label{4.2} \\
                 &\bar{x}_{j}^{\bar{\omega}} \leq q_{j}^{p}\sum_{l \in L} \gamma_{jl}^{p,\bar{\omega}}\bar{y}_{jl} \quad \forall j \in J, \bar{\omega} \in \bar{\Omega} \label{4.3} \\
                 &\bar{x}_{j}^{\bar{\omega}} = \sum_{i \in I} \bar{u}_{ij}^{\bar{\omega}} \quad \forall j \in J, \bar{\omega} \in \bar{\Omega} \label{4.4} \\            
                 &\sum_{j \in J}\bar{x}_{j}^{\bar{\omega}} \leq \sum_{k \in K}d_{k}^{\bar{\omega}} \quad \forall \bar{\omega} \in \bar{\Omega} \label{4.5} \\
                 &\sum_{k \in K} \bar{v}_{jk}^{\bar{\omega}} \leq q_{j}^{p}\sum_{l \in L}\gamma_{jl}^{p,\bar{\omega}} \bar{y}_{jl} \quad \forall j \in J, \bar{\omega} \in \bar{\Omega} \label{4.6}\\
                 & \sum_{k \in K}\bar{v}_{jk}^{\bar{\omega}} + \bar{e}_j^{\bar{\omega}} + \sum_{j^{\prime} \neq j \in J}\bar{\tau}_{jj^{\prime}}^{\bar{\omega}} = \sum_{i \in I}\bar{u}_{ij}^{\bar{\omega}}+\sum_{j^{\prime} \neq j \in J}\bar{\tau}_{j^{\prime}j}^{\bar{\omega}} \quad \forall j \in J, \bar{\omega} \in \bar{\Omega} \label{4.7} \\
                 & \sum_{j^{\prime}\neq j \in J} \bar{\tau}_{j^{\prime}j}^{\bar{\omega}} \leq M_3 \sum_{l \in L}\bar{y}_{jl} \quad \forall j \in J, \bar{\omega} \in \bar{\Omega}\label{4.8}\\
                 & \sum_{j \in J} \bar{v}_{jk}^{\bar{\omega}} \leq d_k^{\bar{\omega}} \quad \forall k \in K, \bar{\omega} \in \bar{\Omega} \label{4.9}\\
                 & \bar{y}_{jl} \in \{0,1\} \quad \forall j \in J, l \in L \label{4.10} \\
                 & \bar{u}_{ij}^{\bar{\omega}}, \bar{v}_{jk}^{\bar{\omega}}, \bar{x}_j^{\bar{\omega}}, \bar{e}_j^{\bar{\omega}}, \bar{\tau}_{jj^{\prime}}^{\bar{\omega}} \geq 0 \quad \forall i \in I, j \in J, j^\prime \neq j \in J, k \in K, \bar{\omega} \in \bar{\Omega} \label{4.11}
            \end{align} 
\end{subequations}    
The objective function \eqref{OF4} maximizes the expected profits when there are no disruptions. Constraints \eqref{4.1} ensure that at most one plant is open in every country. Constraints \eqref{4.2}, \eqref{4.3} and \eqref{4.6} are capacity constraints, limiting the purchases of raw materials, drug production plan and final drug production, respectively. Constraints \eqref{4.5} and \eqref{4.9} are demand constraints. The balance constraints are \eqref{4.4} and \eqref{4.7}, balancing production plan with purchases of raw materials and final drug production with raw materials available in each manufacturing plant, respectively. Raw materials transshipment is only allowed between open plants \eqref{4.8}. Finally, constraints \eqref{4.10} and \eqref{4.11} enforce the domain.  

After the model is solved, the nominal global drug availability, $\bar{g}$, is calculated by equation \eqref{eq:g_bar}. This represents the sum of the company's expected production and the exogenous drug production of other companies when there are not disruptions. 
\begin{equation}\label{eq:g_bar}
    \bar{g}=\sum_{\bar{\omega} \in \bar{\Omega}}\varrho^{\bar{\omega}}\sum_{j \in J}\sum_{k \in K}\bar{v}_{jk}^{\bar{\omega}}+\sum_{k \in K} \eta_{k}
\end{equation}

\section{Further Detail on Export Ban Constraints}\label{supp:sec:EBcon_detail}
For clarity, in this section, we walk through how Inequality \eqref{eq:EB_x} is enforced in Model \eqref{1-2stage}. There are two cases.

When the left-hand side of constraint \eqref{2.8-2stage} is \textit{positive}, the projected global drug supply in a scenario $\omega \in \Omega$ is below a threshold $r$ of the nominal global production $\bar{g}$. This relationship forces the binary auxiliary variable $z^\omega$ to be zero, i.e., countries in set $\bar{J}$ impose export bans (modeled via constraints \eqref{2.9-2stage}).

Conversely, when the left-hand side of constraint \eqref{2.8-2stage} is \textit{non-positive}, there is not a major projected shortage in the global drug supply. In this case, the binary auxiliary variable $z^\omega$ can take the value of zero or one.
\begin{itemize}[noitemsep,topsep=0pt]
    \item If there exists $j^\prime \in \bar{J}$ such that $y_{j^\prime l}=1$ for any $l \in L$ and there is economically attractive demand to meet, then $z^{\omega} = 1$ by constraints \eqref{2.9-2stage} and the sense of the objective function \eqref{OF2-2stage}, maximization of profits.
    \item If $\nexists j^\prime \in \bar{J}$ such that $y_{j^\prime l}=1$ for any $l \in L$, then $z^{\omega}$ can be zero or one. In this case, the value that $z^\omega$ takes will not affect model results since it will not affect the objective function and the international flows. 
\end{itemize}

\section{Data}\label{supp:Data}
This section presents a description of each dataset used for the sets and parameters.

\subsection{\emph{Sets}}

The sets of suppliers $(I)$, potential locations of manufacturing plants $(J)$, and demand countries $(K)$ are from \citeSupp{s-SABOGAL2025}. The set of manufacturing quality levels of facilities $(L)$ is defined as $\{1,2,3\}$, where $l=3$ represents the highest level of quality. Set $\bar{J}$, following procedures presented by \citeSupp{s-SABOGAL2025}, corresponds to countries that have banned the export of drugs according to the World Trade Organization databases \citepSupp{s-WTO2022a,s-WTO2022b}.

\subsection{\emph{Stochastic Parameters}}
Supplier capacity strains $(\gamma_{i}^{s})$ and demand $(d_{k})$ are obtained from \citeSupp{s-SABOGAL2025}. The production capacity strains $(\gamma_{jl}^{p})$ are estimated as follows. For quality level $l=2$, $\gamma_{j2}^{p}$ is set equal to the probability mass functions estimated by \citeSupp{s-SABOGAL2025}. For $l=1$ and $l=3$, probability mass functions of $\gamma_{j2}^{p}$ are modified such that $\mathbb{E}[\gamma_{j2}^{p,\omega}]$ is decreased by 1\% and increased by 1\%, respectively.

\textit{Disruption risks:} Supplier quality disruptions $(\xi_{i}^{qs})$ are modeled with Bernoulli distributions as presented by \citeSupp{s-SABOGAL2025}. Manufacturing plant quality disruptions $(\xi_{jl}^{qp})$ are estimated as follows. For quality level $l=2$, $\xi_{j2}^{qp}$ are set equal to Bernoulli distributions estimated by \citeSupp{s-SABOGAL2025}. For $l=1$ and $l=3$, $\xi_{j1}^{qp} \sim \text{Bernoulli} (\rho_{j1}^{qp})$, and $\xi_{j3}^{qp} \sim \text{Bernoulli}(\rho_{j3}^{qp})$, where $\rho_{j1}^{qp}$ and $\rho_{j3}^{qp}$ are estimated using the disruptions per year rate presented by \citeSupp{s-Tucker2020} increased by 90\% and decreased by 90\%, respectively. 

Data to parameterize the natural disaster disruptions ($\xi_k^{nd}$) approach from Section (S\ref{supp:stochastic})
are given as follows.Level (1): the occurrence of a natural disaster in country $k^\prime \in K$, is modeled with the Bernoulli random variable $\mathcal{U}_{k^\prime}$, that equals zero with probability $\psi_{k^\prime}$ indicating that country $k^\prime$ is the origin of a natural disaster. $\psi_{k^\prime}$ is estimated using The International Disaster Database EM-DAT \citepSupp{s-EM-DAT2024}. We obtained information on natural disasters that occurred around the world between the years 2000-2022, and that affected large populations (more than 500000 people). We considered the top hazards that affect pharmaceutical supply chains, i.e., flood, storm, wildfire and earthquake \citepSupp{s-pharmamanufacturing}.

Level (2): the propagation effects of natural disasters exist between countries with common hazards. The set $H_{k}$ for every $k \in K$, represents the set of countries with common hazards with country $k$, and is defined based on proximity. At the base case, for country $k$ we consider the countries that are at most 3218 km (2000 miles) from country $k$. According to the EM-DAT database \citepSupp{s-EM-DAT2024}, most natural disasters that impact multiple countries occur between countries that are less than 2000 miles apart. 

We set the baseline strength of the disaster propagation to be $\zeta= 1000$ km as base case (medium propagation), which is in the same order of magnitude as other works in this area (e.g., \citeSupp{s-Lu-Ran2015,s-Li2018}). We evaluate the effects of a high propagation value ($\zeta = 1300$ km) in the Set III instances.

Level (3): When a country is impacted by a natural disaster, it is not always the case that the production facilities will be disrupted. We consider the probability that manufacturing capacity in country $k$ is disrupted given that a natural disaster hits country $k$. It is parameterized according to the same set of natural disasters from Level (1). We estimate the conditional probability as the average of the affected population in country $k$ vs. the total population of the country $k$ across all included disasters that affected country $k$ \citepSupp{s-EM-DAT2024}.

\subsection{\emph{Deterministic Parameters}}

Raw materials costs ($c^{s}$), production costs ($c^{p}$), transportation costs ($c^{t}$), prices $(p)$, supplier capacity $(q^{s})$, manufacturing plants capacity $(q^{p})$, and exogenous drug exports $(\eta)$ are taken from \citeSupp{s-SABOGAL2025}.

Fixed costs for opening and operating manufacturing plants ($c_{jl}^{f})$ are estimated as follows. For quality level $l=2$, $c_{j2}^{f}$ are set equal to the data provided by \citeSupp{s-SABOGAL2025}. For $l=1$, fixed costs are set equal to $c_{j2}^{f}$ reduced by 25\% of the change in quality with respect to $l=2$ (-90\%). For $l=3$, fixed costs are set equal to $c_{j2}^{f}$ increased by 50\% of the improvement in quality with respect to $l=2$ (90\%). The holding costs of raw materials inventory ($c^h$) are set as 36\% of the price \citepSupp{s-Tucker2020}.

One value for the threshold to impose export bans ($r$) is set to 0.65 based on \citepSupp{s-bulgaria}. In Bulgaria, if the availability of a specific drug is less than 65\% of the domestic demand for a period of one month, local authorities may implement export bans on that drug. Additional cases $r\in\{0.8,0.95\}$ are included to evaluate the effects of higher sensitivity.

\section{Additional Details of Algorithms}\label{supp:Algorithms}

\subsection{\emph{Alternating Disjunctive Cut Method}}
A cut-generating linear program \eqref{Separation_problem} that reflects the dynamics presented in Section \ref{sub:disjuctive}, together with a normalizing constraint \eqref{normalizing-const}, is used to find the multipliers and coefficients of the disjunctive cut.
\begin{subequations}\label{Separation_problem}
    \small
        \begin{align}
             & \max \quad \hat{\theta}^{\omega} \bar{\delta}^{\omega} + \sum_{j \in J} \sum_{l \in L} \hat{y}_{jl}\bar{\alpha}_{jl}^{\omega} - \bar{\beta}^{\omega}\\
             \text{s.t.  }
             & \beta_{0}^{\omega}\lambda_{0}^c + \sum_{j \in J} \lambda_{0j}^{m} + \sum_{j \in J}\sum_{l \in L} \lambda_{0jl}^u \leq \bar{\beta}^{\omega} \\
             & \beta_{1}^{\omega}\lambda_{1}^c + \sum_{j \in J} \lambda_{1j}^{m} + \sum_{j \in J}\sum_{l \in L} \lambda_{1jl}^u \leq \bar{\beta}^{\omega} \\
             & \bar{\alpha}_{jl}^{\omega} \leq -\alpha_{0jl}^{\omega}\lambda_{0}^c  + \lambda_{0j}^{m} + \lambda_{0jl}^{u} \quad \forall j \in J, l \in L \\
             & \bar{\alpha}_{jl}^{\omega} \leq -\alpha_{1jl}^{\omega}\lambda_{1}^c  + \lambda_{1j}^{m} + \lambda_{1jl}^{u} \quad \forall j \in J, l \in L \\
             & \bar{\delta}^{\omega} \leq \lambda_{0}^c \\
             & \bar{\delta}^{\omega} \leq \lambda_{1}^c \\
             & \lambda_{0}^c + \lambda_{1 }^c = 1 \label{normalizing-const} \\
             & \lambda_{0}^{c}, \lambda_{1}^{c}, \lambda_{0j}^{m}, \lambda_{1j}^{m}, \lambda_{0jl}^{u}, \lambda_{1jl}^{u} \geq 0 \quad \forall j \in J, l \in L
        \end{align} 
\end{subequations}

Algorithm \ref{Algorithm:D} presents the Alternating Disjunctive Cut Method that alternates between Benders optimality cuts and disjunctive cuts to solve model \eqref{1-2stage}.

\begin{algorithm}[h!]
\caption{Alternating Disjunctive Cut Method}\label{Algorithm:D}
\begin{algorithmic}[1]
\spacingset{1}
\small
\STATE \textbf{Input:} Incumbent solution $(\boldsymbol{\hat{\theta}}, \mathbf{\hat{y}})$ of $(MP)$
\FOR{$\omega \in \Omega$}
    \STATE Solve $\mathbf{Q_{LP}^{\omega}(\hat{y})}$ 
    \IF{$\hat{\theta}^{\omega} > \mathbf{Q_{LP}^{\omega}(\hat{y})}$}
        \STATE Add a Benders optimality cut \eqref{benderscut} to set $\mathcal{O}$ 
    \ELSE
        \STATE Solve $\mathbf{Q_{LP}^\omega}(\mathbf{\hat{y}}, z^\omega=0)$
        \STATE Calculate Benders optimality cut coefficients $(\alpha^{\omega}_{0},\beta^{\omega}_{0})$ with \eqref{eq:B0}-\eqref{eq:alphaz}\\
        \STATE \small Solve $\mathbf{Q_{LP}^\omega}(\mathbf{\hat{y}},z^\omega=1)$
        \IF{$\mathbf{Q_{LP}^\omega}(\mathbf{\hat{y}},z^\omega=1)$ is feasible}
            \STATE Calculate Benders optimality cut coefficients $(\alpha^{\omega}_{1},\beta^{\omega}_{1})$ with \eqref{eq:B0}-\eqref{eq:alphaz}
            \STATE $\mathbf{Q^\omega(\hat{y}}) \leftarrow \max\{\mathbf{Q_{LP}^\omega}(\mathbf{\hat{y}}, z^\omega=0),\mathbf{Q_{LP}^\omega}(\mathbf{\hat{y}}, z^\omega=1)\}$ 
            \STATE Solve cut-generating linear program $\eqref{Separation_problem}$ using ($\hat{\theta}^{\omega}, \mathbf{\hat{y}},\alpha^{\omega}_{0},\beta^{\omega}_{0}, \alpha^{\omega}_{1},\beta^{\omega}_{1})$ to find ($\bar{\delta}^{\omega}$,$\bar{\alpha}^{\omega}$,$\bar{\beta}^{\omega}$) coefficients of disjunctive cut \eqref{eq:disjuctive-cut-w}
            \STATE $\alpha^{\omega} \leftarrow -\bar{\alpha}^{\omega}/\bar{\delta}^{\omega}$, $\beta^{\omega} \leftarrow \bar{\beta}^{\omega}/\bar{\delta}^{\omega}$  
        \ELSE
            \STATE $\mathbf{Q^\omega(\hat{y})} \leftarrow \mathbf{Q_{LP}^\omega}(\mathbf{\hat{y}}, z^\omega=0), \alpha^{\omega} \leftarrow \alpha^{\omega}_0, \beta^{\omega}  \leftarrow \beta^{\omega}_0$ 
        \ENDIF
        \IF{$\hat{\theta}^{\omega} > \mathbf{Q^{\omega}(\hat{y})}$} 
            \STATE Add cut $\theta^\omega \leq \alpha^\omega \mathbf{y}+\beta^\omega$ to set $\mathcal{O}$ 
         \ENDIF
    \ENDIF
\ENDFOR
\end{algorithmic}
\end{algorithm}
\newpage

\subsection{\emph{Alternating Bilinear Cut Method}}
\begin{subequations}\label{eq:new_mp_constraints}

As discussed in Section \ref{sub:BilinearCuts}, constraints \eqref{lin_1}--\eqref{lin_4} are used to linearize $y_{jl}z^\omega$ (a product of binary variables) through the use of an auxiliary binary variable, $w_{jl}^\omega$.  
\small
    \begin{align}
        & w_{jl}^{\omega} \leq z^{\omega} \quad \forall j \in J, l \in L, \omega \in \Omega \label{lin_1} \\
        & w_{jl}^{\omega} \leq y_{jl} \quad \forall j \in J, l \in L, \omega \in \Omega \label{lin_2} \\
        & w_{jl}^{\omega} \geq z^{\omega} + y_{jl} - 1 \quad \forall j \in J, l \in L, \omega \in \Omega \label{lin_3} \\
        & w_{jl}^{\omega} \geq 0 \quad \forall j \in J, l \in L, \omega \in \Omega \label{lin_4}
    \end{align}
\end{subequations}

Algorithm \ref{Algorithm:B} is the Alternating Bilinear Cut Method. It alternates between Benders optimality cuts and linearized bilinear cuts.
\begin{algorithm}[h!]
\caption{Alternating Bilinear Cut Method}\label{Algorithm:B}
\begin{algorithmic}[1]
\spacingset{1}
\small
\STATE \textbf{Input:} Incumbent solution $(\boldsymbol{\hat{\theta}}, \mathbf{\hat{y}})$ of $(MP^\prime)$
\FOR{$\omega \in \Omega$}
    \STATE Solve $\mathbf{Q_{LP}^{\omega}(\hat{y})}$ 
    \IF{$\hat{\theta}^{\omega} > \mathbf{Q_{LP}^{\omega}(\hat{y})}$}
        \STATE Add a Benders optimality cut \eqref{benderscut} to set $\mathcal{O^\prime}$ 
    \ELSE
        \STATE Solve $\mathbf{Q_{LP}^\omega}(\mathbf{\hat{y}}, z^\omega=0)$
        \STATE Calculate Benders optimality cut coefficients $(\alpha^{\omega}_{0},\beta^{\omega}_{0})$ with \eqref{eq:B0}-\eqref{eq:alphaz}\\
        \STATE \small Solve $\mathbf{Q_{LP}^\omega}(\mathbf{\hat{y}},z^\omega=1)$
        \IF{$\mathbf{Q_{LP}^\omega}(\mathbf{\hat{y}},z^\omega=1)$ is feasible}
            \STATE Calculate Benders optimality cut coefficients $(\alpha^{\omega}_{1},\beta^{\omega}_{1})$ with \eqref{eq:B0}-\eqref{eq:alphaz}
            \STATE $\mathbf{Q^\omega(\hat{y}}) \leftarrow \max\{\mathbf{Q_{LP}^\omega}(\mathbf{\hat{y}}, z^\omega=0),\mathbf{Q_{LP}^\omega}(\mathbf{\hat{y}}, z^\omega=1)\}$ 
        \ELSE
            \STATE $\mathbf{Q^\omega(\hat{y})} \leftarrow \mathbf{Q_{LP}^\omega}(\mathbf{\hat{y}}, z^\omega=0), \alpha_{1}^{\omega} \leftarrow 0, \beta_{1}^{\omega}  \leftarrow 0$
        \ENDIF
        \IF{$\hat{\theta}^{\omega} > \mathbf{Q^{\omega}(\hat{y})}$} 
            \STATE Add cut \eqref{master_bilinearcut} to set $\mathcal{O^\prime}$ 
         \ENDIF
    \ENDIF
\ENDFOR
\end{algorithmic}
\end{algorithm}

\newpage

\section{Additional Computational Results}\label{supp:comp_results}

\subsection{\emph{Single-Cut Adaptations of Our Cutting Plane Methods}}
In this subsection, we present the computational performance of the single-cut variants of the proposed algorithms. In these versions, a single cut is intended to be added to the master problem every time a feasible binary master problem solution is found. Table \ref{supp:table_singlecut} summarizes the computational results for these modified methods. For each method, the table reports the average solution time (seconds) and the average number of branch-and-bound nodes explored of the samples solved to optimality within the time limit. These are under columns “Time” and “Nodes”, respectively. The column “Solved” indicates the number of samples (out of five) that were solved to optimality within the time limit. The column “Gap” reports the average gap of the samples that were not solved to optimality.

\begin{table}[ht!]
\renewcommand{\arraystretch}{1.1}
\scriptsize
\centering
\setlength{\tabcolsep}{4pt} 
\caption{Computational performance for Set I (baseline) across single-cut methods}
\begin{tabular}{c|c|cccc|cccc|cccc}
    \hline
    \multicolumn{2}{c|} {\textbf{Set I}} & \multicolumn{4}{c|}{\textbf{Alt. Integer L-shaped}} & \multicolumn{4}{c|}{\textbf{Alt. Disjunctive}} & \multicolumn{4}{c}{\textbf{Alt. Bilinear}} \\
    \hline
    $r$ & \textbf{$|\Omega|$} & \textbf{Time} & \textbf{Nodes} & \textbf{Solved} & \textbf{Gap} & \textbf{Time} & \textbf{Nodes} & \textbf{Solved} & \textbf{Gap} & \textbf{Time} & \textbf{Nodes} & \textbf{Solved} & \textbf{Gap} \\
    \hline
        \multirow{3}{*}{0.65} & 100 & 101 & 1262 & 5 & 0.0\% & 115 & 1086 & 5 & 0.0\% & 120 & 1196 & 5 & 0.0\% \\
        & 500 & 695 & 1228 & 5 & 0.0\% & 742 & 1188 & 5 & 0.0\% & 793 & 1010 & 5 & 0.0\% \\
        & 1000 & 1588 & 1763 & 5 & 0.0\% & 1923 & 1438 & 5 & 0.0\% & 2000 & 1389 & 5 & 0.0\% \\
        \hdashline 
        \multirow{3}{*}{0.80} & 100 & 102 & 1262 & 5 & 0.0\% & 111 & 1086 & 5 & 0.0\% & 122 & 1196 & 5 & 0.0\% \\
        & 500 & 677 & 1228 & 5 & 0.0\% & 744 & 1188 & 5 & 0.0\% & 771 & 1010 & 5 & 0.0\% \\
        & 1000 & 1593 & 1763 & 5 & 0.0\% & 1841 & 1438 & 5 & 0.0\% & 1961 & 1389 & 5 & 0.0\% \\
        \hdashline 
        \multirow{3}{*}{0.95} & 100 & 159 & 1071 & 5 & 0.0\% & 198 & 1191 & 5 & 0.0\% & 249 & 1547 & 5 & 0.0\% \\
        & 500 & 878 & 1457 & 5 & 0.0\% & 1007 & 1356 & 5 & 0.0\% & 2541 & 2986 & 5 & 0.0\% \\
        & 1000 & 2446 & 2029 & 4 & 0.8\% & 2841 & 2085 & 4 & 0.2\% & - & - & 0 & 3.1\% \\
    \hline
\end{tabular}
\label{supp:table_singlecut}
\end{table}

\subsection{\emph{Alternating Integer L-Shaped Method with Disjunctive and Bilinear Cuts}}

We implemented two additional methods: (i) a combination of Integer L-shaped cuts with Disjunctive cuts and (ii) a combination of Integer L-shaped cuts with Bilinear cuts. Both methods were evaluated using instances in Set I with a geopolitical threshold of $r=0.95$. For these experiments, we tested the multi-cut versions of the proposed approaches and used the alternating strategy with the Benders optimality cuts. Computational results are presented in Table \ref{table:AILDC-BC}.

\begin{table}[ht!]
\centering
\footnotesize
\caption{Computational performance of combined methods for Set I (baseline) with $r=0.95$}
\renewcommand{\arraystretch}{1.1}
\setlength{\tabcolsep}{8pt} 
\begin{tabular}{c|ccc|ccc}
\hline
{\textbf{Set I}} & \multicolumn{3}{c|}{\textbf{Alt. Integer L-shaped + Disjunctive cuts}} & \multicolumn{3}{c}{\textbf{Alt. Integer L-shaped + Bilinear cuts}} \\
\hline
\textbf{$|\Omega|$} & \textbf{Time} & \textbf{Nodes} & \textbf{Solved} & \textbf{Time} & \textbf{Nodes} & \textbf{Solved}\\
\hline
100 & 139 & 365 & 5 & 271 & 609 & 5 \\
500 & 792 & 497 & 5 & 1168 & 719 & 5 \\
1000 & 1761 & 501 & 5 & 2239 & 617 & 5 \\
\hline
\end{tabular}
\label{table:AILDC-BC}
\end{table}

\subsection{\emph{Additional Results of the Alternating Integer L-Shaped Method}}

In this subsection, we present computational results for variants of the Alternating Integer L-shaped Cut Method. As explained in the manuscript, when variable $z^\omega$ is binary in the optimal solution of $\mathbf{Q_{LP}^{\omega}}\mathbf{(\hat{y})}$, then an integer L-shaped cut can be generated as a byproduct when attempting to add a Benders optimality cut for scenario $\omega \in \Omega$. Considering this situation, we tested adding both cuts (Benders and Integer L-shaped), adding only a Benders cut, and adding only an Integer L-shaped cut. Table \ref{tab:AIL_variants} reports the results of these experiments. We used three independent samples of 1000 scenarios under Set I and $r = 0.95$. All experiments were conducted with a time limit of 3600 seconds.

\begin{table}[h]
\centering
\footnotesize
\caption{Computational results of Alternating Integer L-shaped Method variants with $|\Omega|=1000$}
\renewcommand{\arraystretch}{1.1}
\begin{tabular}{c|cc|cc|cc}
\hline
{$r=0.95$} & \multicolumn{2}{c|}{\textbf{Both Cuts}} & \multicolumn{2}{c|}{\textbf{Only Benders}} & \multicolumn{2}{c}{\textbf{Only Integer L-shaped}} \\
\hline
\textbf{Sample} & \textbf{Time} & \textbf{Gap} & \textbf{Time} & \textbf{Gap} & \textbf{Time} & \textbf{Gap} \\
\hline
a & 1154 & 0.0\% & 1183 & 0.0\% & $>$3600 & 1.6\%  \\
b & 1081 & 0.0\% & 1314 & 0.0\% & 2188 & 0.0\%     \\
c & 1089 & 0.0\% & 1153 & 0.0\% & $>$3600 & 2.4\%  \\
\hline
\end{tabular}
\label{tab:AIL_variants}
\end{table}

\subsection{\emph{Additional Results of the Alternating Bilinear Cut Method}}

In this subsection, we present computational results for a variant of the Alternating Bilinear Cut Method. In this variant, we include a feasibility cut for scenario $\omega \in \Omega$ when solving $\mathbf{Q_{LP}^\omega}(\mathbf{\hat{y}}, z^\omega=1)$ is infeasible. The feasibility cut that is added to the master problem $(MP^\prime)$ is presented in \eqref{eq:feasibility_cut}. Table \ref{tab:AB_variant} reports the results of this variant and compares them with the original implementation.  
\begin{equation}\label{eq:feasibility_cut}
     z^{\omega} \leq |\bar{Y}_{1}| - \sum_{(j,l) \in \bar{Y}_{1}} y_{jl} + \sum_{(j,l) \in \bar{Y}_{0}} y_{jl} 
\end{equation}
Where, $\bar{Y}_{1}=\{(j,l): y_{jl} =1$ for all $j \in J, l \in L\}$ and $\bar{Y}_{0}=\{(j,l): y_{jl} =0$ for all $j \in J, l \in L\}$

\begin{table}[h!]
\centering
\footnotesize
\renewcommand{\arraystretch}{1}
\caption{Computational results for the Alternating Bilinear Cut Method variants under Set I}
\begin{tabular}{c|c|ccc|ccc}
\hline
\multicolumn{2}{c|}{\textbf{Alt. Bilinear Cut}} & \multicolumn{3}{c|}{\textbf{Without Feasibility Cut}} & \multicolumn{3}{c}{\textbf{With Feasibility Cut}} \\
\hline
{$r$} & {$|\Omega|$}  & \textbf{Time} & \textbf{Nodes} & \textbf{Solved} & \textbf{Time} & \textbf{Nodes} & \textbf{Solved} \\
\hline
\multirow{3}{*}{0.65} & 100  & 245  & 612 & 5 & 242  & 612 & 5 \\
     & 500  & 928  & 666 & 5 & 972  & 666 & 5 \\
     & 1000 & 1387 & 567 & 5 & 1370 & 567 & 5 \\
\hdashline
\multirow{3}{*}{0.80}  & 100  & 255  & 607 & 5 & 247  & 607 & 5 \\
     & 500  & 837  & 715 & 5 & 952  & 659 & 5 \\
     & 1000 & 1493 & 568 & 5 & 1483 & 512 & 5 \\
\hdashline
\multirow{3}{*}{0.95} & 100  & 302  & 633 & 5 & 279  & 617 & 5 \\
     & 500  & 1116 & 630 & 5 & 1194 & 667 & 5 \\
     & 1000 & 1935 & 632 & 5 & 2507 & 595 & 5 \\
\hline
\end{tabular}
\label{tab:AB_variant}
\end{table}

\newpage